\documentclass[11pt]{article}
\usepackage[noadjust]{cite}

\usepackage[margin=1in]{geometry}
\usepackage{setspace}
\usepackage{color}

\usepackage{amsmath}   
\usepackage{amssymb}   
\usepackage{amsthm}     
\usepackage{mathtools} 
\usepackage{dsfont}  
\usepackage{accents}

\usepackage{graphicx}
\usepackage{epstopdf}
\usepackage{booktabs}  
\usepackage{adjustbox}
\usepackage{multirow}  
\usepackage{caption}
\usepackage{subcaption}

\usepackage{tikz}
\usetikzlibrary{shapes,arrows}
\usepackage{tikz-cd}

\usepackage{enumerate}
\usepackage[ruled,vlined]{algorithm2e}

\usepackage{authblk}    
\usepackage[normalem]{ulem}  
\usepackage{textcomp}  
\usepackage{appendix}

\usepackage[colorlinks=true, allcolors=blue]{hyperref}

\usepackage{float}
\makeatletter
\newcases{centercases}{\quad}
  {\hfil\m@th\displaystyle{##}\m@th\displaystyle{##}\hfil}
  {\m@th\displaystyle{##}\m@th\displaystyle{##}\hfil}{\lbrace}{.}
\makeatother

\usepackage{grffile}
\usetikzlibrary{decorations.pathreplacing}
\usepackage{wrapfig}
\usepackage{csquotes}

\newcommand{\commentout}[1]{}

\theoremstyle{plain}

\theoremstyle{definition}

\newlength{\dhatheight}

\begin{document}
	
	\title{PriorIDENT: Prior-Informed PDE Identification from Noisy Data}
	\author{
 Cheng Tang\thanks{Department of Mathematics, Hong Kong Baptist University, Kowloon Tong, Hong Kong. Email: 22481184@life.hkbu.edu.hk. }, 
 Hao Liu\thanks{Department of Mathematics, Hong Kong Baptist University, Kowloon Tong, Hong Kong.
		Email: haoliu@hkbu.edu.hk.},  
 Dong Wang\thanks{School of Science and Engineering, The Chinese University of Hong Kong (Shenzhen), Shenzhen, Guangdong 518172, China \& Shenzhen International Center for Industrial and Applied Mathematics, Shenzhen Research Institute of Big Data, Guangdong 518172, China \& Shenzhen Loop Area Institute, Guangdong 518048, China. Email: wangdong@cuhk.edu.cn}
  }
	
	\date{}
	\maketitle

	\begin{abstract}
Identifying governing partial differential equations (PDEs) from noisy spatiotemporal data remains challenging due to differentiation-induced noise amplification and ambiguity from overcomplete libraries. We propose a prior-informed weak-form sparse-regression framework that resolves both issues by refining the dictionary before regression and shifting derivatives onto smooth test functions. Our design encodes three compact physics priors—Hamiltonian (skew-gradient and energy-conserving), conservation-law (flux-form with shared cross-directional coefficients), and energy-minimization (variational, dissipative)—so that all candidate features are physically admissible by construction. These prior-consistent libraries are coupled with a subspace-pursuit pipeline enhanced by trimming and residual-reduction model selection to yield parsimonious, interpretable models. Across canonical systems—including Hamiltonian oscillators and the three-body problem, viscous Burgers and two-dimensional shallow-water equations, and diffusion and Allen--Cahn dynamics—our method achieves higher true-positive rates, stable coefficient recovery, and structure-preserving dynamics under substantial noise, consistently outperforming no-prior baselines in both strong- and weak-form settings. The results demonstrate that compact structural priors, when combined with weak formulations, provide a robust and unified route to physically faithful PDE identification from noisy data.
\end{abstract}

\section{Introduction}
The discovery of governing physical laws from data has attracted increasing attention in recent years, as it bridges scientific modeling and data-driven learning. In particular, learning partial differential equations (PDEs) directly from spatio-temporal observations offers new opportunities for understanding complex systems in physics, engineering, and biology. A large body of work—ranging from sparse-regression~\cite{Discoveringgoverning_sparsereg,Schaeffer2017LearningPD_sparsereg,Data_drivenPDE_sparsereg,IDENT_sparsereg,he2022robust, Messenger_2021_sparse_reg, Rudy_sparse_reg, tang2022weakident,  WU2019200, cheng2025identification, Messenger_2021}  to neural PDE/operator approaches~\cite{long2018pdenet,long2019pde, NEURIPS2023_70518ea4,chen2022deep,churchill2025principal,wu2020data}—has advanced the field of equation discovery. More recently, the scope of these methods has been expanded to address more complex dynamics and representations: Stoch-IDENT~\cite{cui2025stoch} was proposed for identifying stochastic PDEs, Phase-IDENT~\cite{yang2026phase} addressed dynamical systems with phase transitions, and spectral-based frameworks were developed for identification in the Fourier domain~\cite{tang2026identifying, heinrich2025rediscovering}. A review of related methods can be found in \cite{he2025recent}. Regression-based approaches typically assume a generic form
\begin{equation} \label{pde_form}
\partial_t u = F(\textbf{x}, \partial_\textbf{x} u, \partial^2_\textbf{x} u,\dots),
\end{equation}
where the right-hand side is represented as a combination of terms drawn from a predefined overcomplete feature dictionary containing spatial derivatives and nonlinear interactions. The governing operator \(F\) is modeled as a sparse linear combination of these candidate features, reflecting the empirical observation that physically meaningful PDEs are parsimonious.  
Such sparsity not only enhances interpretability but also mitigates overfitting, making the identification of the active terms and their coefficients a central objective in PDE discovery.

Nevertheless, existing methods are frequently hampered by noise amplification—an artifact of numerical differentiation—and model ambiguity due to overcomplete dictionaries. The redundancy of candidate terms complicates the selection of correct governing components, often resulting in the misselection of non-physical terms that merely fit noise. For instance, when a system evolves by minimizing an energy functional, noise can lead to the discovery of PDE terms that are not gradients of any energy, thereby contradicting fundamental physical principles. This paper therefore focuses on two core challenges: improving robustness against noise and developing a principled framework for dictionary construction to ensure physical interpretability. A promising pathway, explored across various disciplines, is the integration of prior knowledge as a regularizing constraint, which offers a natural foundation for enhancing both the robustness and interpretability of PDE discovery.

In classical image denoising, one formulates an optimization problem
\[
\hat{x} = \arg\min_x \frac{1}{2}\|y - x\|_2^2 + \lambda R(x),
\]
balancing fidelity to the noisy observation \(y\) with a regularization term \(R(x)\) that encapsulates prior assumptions~\cite{fan2019brief}.  
The regularizer \(R(x)\) can be interpreted as the negative log-prior on the true image, and choosing an appropriate prior is therefore crucial. Early image-denoising studies underscore the value of priors: classical smoothness and total-variation regularizations impose structural constraints that suppress noise while preserving essential features.  
Such regularization-based strategies bias solutions toward physically or statistically plausible patterns, a principle that naturally extends to the discovery of governing equations.

A similar philosophy has recently emerged in the context of PDE identification, where physical priors are embedded directly into neural network architectures.  
For example, Lee et al.~\cite{lee2022structure} enforced a Hamiltonian structure within neural ODEs to preserve energy and symplecticity by constraining the learned dynamics to
\[
\dot{x} = J \nabla H(x),
\]
where \(J\) is skew-symmetric.  
Li et al.~\cite{li2023learning} incorporated a conservation-law prior via neural finite-volume schemes, while Eidnes and Lye~\cite{eidnes2024pseudo} introduced pseudo-Hamiltonian networks that disentangle conservative and dissipative effects.  
These studies collectively illustrate that physics-informed priors can greatly improve robustness and interpretability by aligning learned models with underlying physical structures.

Motivated by this perspective, we incorporate prior knowledge directly into the sparse-regression framework for PDE discovery. Instead of embedding priors within neural architectures, we use them to guide dictionary construction, ensuring candidate terms are both physically admissible and interpretable.

In addition to prior-informed modeling, weak-form formulations have proven highly effective for mitigating noise sensitivity.  
By transferring differentiation from discrete data to smooth test functions, the weak form markedly enhances stability and numerical robustness.  
Representative examples include WeakSINDy~\cite{Messenger_2021}, WeakIdent~\cite{tang2022weakident}, and WG-IDENT~\cite{TANG2026114454}, the latter also accommodating spatially varying coefficients through group-structured sparse selection.

Motivated by these insights, we propose a unified \emph{prior-informed weak-form PDE identification framework} that addresses the two central challenges: 
(i) the incorporation of physical priors into sparse-regression discovery, and 
(ii) robustness against noise through weak formulations.  

Our contributions are summarized as follows:
\begin{enumerate}
    \item
    We propose a novel method to incorporate prior information with dictionary design, which improves the overall robustness and ensures that the identified model is physically meaningful.

    \item  
    We develop a unified framework combining prior-informed dictionary construction (Hamiltonian / flux-form / gradient-flow) with a weak-form sparse-regression pipeline, incorporating trimming and residual-based model selection.

    \item
    Extensive experiments show that our method achieves higher successful rates, stable coefficient recovery, and structure-preserving dynamics even under substantial noise.  
    It consistently outperforms existing no-prior baselines in both weak-form~\cite{tang2022weakident} and strong-form~\cite{IDENT_sparsereg} settings.
\end{enumerate}

The rest of this paper is organized as follows:
Section~\ref{sec:pde_identification_framework} presents an overview of the overall PDE identification framework.  
In Section~\ref{sec:prior_dictionary_refinement}, we introduce the proposed prior-informed dictionary refinement strategy, detailing three representative forms of prior knowledge: the Hamiltonian, conservation-law, and energy–dissipation priors.  
We provide illustrative case studies for each prior in Section~\ref{sec:case_studies}, emphasizing the corresponding dictionary construction and weak-form formulations.  
Section~\ref{sec:numericalexp} reports extensive numerical experiments on these systems to evaluate the performance and robustness of the proposed approach.  
We conclude this paper in Section~\ref{sec:conclusions}.

\section{Overview of the PDE identification framework}
\label{sec:pde_identification_framework}
Given measurements of a state variable $\mathbf{u}(x,t) \in \mathbb{R}^d$, where $d$ denotes the dimension of the state vector, the objective is to infer the expression of the underlying PDE:
\[
\mathbf{u}_t = \mathcal{N}(\mathbf{u}, \mathbf{u}_x, \mathbf{u}_{xx}, \dots),
\]
where $\mathcal{N}(\cdot)$ represents an unknown nonlinear differential operator depending on the state and its spatial derivatives. While various approaches exist for PDE identification, here we adopt the framework used in~\cite{he2022robust, tang2022weakident, he2023group, TANG2026114454, tang2026identifying}, which formulates the identification task as an $\ell_0$ minimization problem. The resulting identification procedure comprises five primary steps.

\textbf{Step 1: Dictionary construction.}
First, a candidate library of $P$ potential functional terms, denoted by $\{\Phi_j\}_{j=1}^P$, is constructed. These terms typically consist of polynomial combinations of the state and its spatial derivatives (e.g., $u, u^2, u_x, u_{xx}$). By evaluating each $\Phi_j$ across the $M$ spatiotemporal observations, we assemble the design matrix $\mathbf{\Theta} \in \mathbb{R}^{M \times P}$, where the $j$-th column corresponds to the discrete trajectory of $\Phi_j$. Organizing the temporal derivatives (or their weak-form projections) into a target vector $\mathbf{b} \in \mathbb{R}^M$ casts the identification problem as a linear system:
\[
\mathbf{b} = \mathbf{\Theta}\mathbf{c},
\]
where the unknown vector $\mathbf{c} \in \mathbb{R}^P$ encodes the active physical terms. Since we expect the true PDE is simple, i.e., it relies on only a small subset of the candidate library, the coefficient vector $\mathbf{c}$ is assumed to be sparse.

\textbf{Step 2: Sparse regression via subspace pursuit.}
Based on the sparsity assumption established above, we employ the \emph{Subspace Pursuit (SP)} algorithm~\cite{he2022robust} to recover the coefficient vector $\mathbf{c}$. Given a sparsity $\theta$, SP finds a $\theta$-sparse vector that best fits the data. We use SP to generates a sequence of candidate models with distinct sparsity.

\textbf{Step 3: Trimming of marginal features.}
The support set of a candidate model may include spurious terms with negligible contributions, particularly when the sparsity $\theta$ exceeds the true sparsity $\theta^*$. A \emph{trimming} procedure~\cite{tang2022weakident} is applied to mitigate such over-selection. For a model with sparsity $\theta$, let $\mathbf{\Theta}_v^{\theta}$ denote the column of $\mathbf{\Theta}$ corresponding to the $v$-th feature of this model. We define a contribution score $\chi_v^\theta$ for the $v$-th feature as:
\begin{equation}
\chi_v^\theta = \frac{\alpha_v^\theta}{\max_{1\leq v' \leq \theta} \alpha_{v'}^\theta}, \quad \text{where} \quad \alpha_v^\theta = \|\mathbf{\Theta}_v^\theta \mathbf{c}_v^\theta\|_2, \quad v = 1, \dots, \theta.
\label{eq:contribution_score}
\end{equation}
Here, $\alpha_v^\theta$ quantifies the contribution of the $v$-th term to the approximation of $\mathbf{b}$, while $\chi_v^\theta$ represents its relative counterpart. Features satisfying $\chi_v^\theta < \tau$ for a fixed threshold $\tau > 0$ are discarded. Subsequently, the coefficients of the remaining features are updated via least squares. This procedure yields a pool of refined PDE candidates across all sparsity levels.

\textbf{Step 4: Model selection via Reduction in Residual (RR).}

The optimal sparsity level $\theta^*$ is selected using the \emph{Reduction in Residual (RR)} criterion~\cite{he2023group}. For each candidate sparsity level $\theta$, let $R_\theta = \|\mathbf{\Theta} \mathbf{c}_\theta - \mathbf{b}\|_2^2$ denote the squared residual of the least-squares fit corresponding to the support $\mathcal{S}_\theta$. The average reduction ratio,
\[
s_\theta = \frac{R_\theta - R_{\theta+L}}{L R_1},
\]
quantifies the normalized improvement in accuracy with increasing sparsity, where $L \ge 1$ is a fixed integer. The optimal sparsity level $\theta^*$ is determined as the smallest $\theta$ for which $s_\theta$ falls below a predefined tolerance $\rho^{\mathrm{R}}$, indicating marginal returns in data fitting. This criterion balances parsimony and accuracy.

\textbf{Step 5: Final identification and reconstruction.}

We obtain the final coefficient vector $\mathbf{c}^*$ by least squares regression using partial columns of $\mathbf{\Theta}$ corresponding to the optimal support:
\[
\min_{\mathbf{c}} \|\mathbf{\Theta}\mathbf{c} - \mathbf{b}\|_2^2 \quad \text{subject to } \operatorname{supp}(\mathbf{c}) = \operatorname{supp}(\mathbf{c}_{\theta^*}),
\]
where $\mathbf{c}_{\theta^*}$ denotes the candidate coefficient vector at the selected sparsity level $\theta^*$. The non-zero entries of $\mathbf{c}^*$ determine the identified PDE structure.

\section{Prior-informed dictionary refinement for PDE identification}
\label{sec:prior_dictionary_refinement}

The framework discussed in Section \ref{sec:pde_identification_framework} offers a systematic data-driven approach, whose accuracy and interpretability depend on the quality of the dictionary. Treating $\mathcal{N}$ as unknown necessitates the construction of large, redundant feature libraries to ensure completeness. However, such overparameterization can induce ill-conditioning, amplify measurement noise, and yield physically inconsistent models. 

In this paper, we propose a \emph{prior-informed dictionary refinement} strategy to mitigate these issues by embedding physical knowledge into the candidate space construction. Structural priors, such as Hamiltonian formulations, conservation laws, or energy-dissipation principles, are incorporated to constrain admissible features. This refinement restricts the search space to features consistent with physical laws, enhancing robustness against noise and ensuring fidelity to the underlying dynamics.

\subsection{Proposed framework}
Our framework integrates physical knowledge directly into the identification pipeline by imposing a set of structural constraints $\mathcal{P}$ (e.g., symmetries, conservation laws, or Hamiltonian/gradient-flow structure) on the governing operator $\mathcal{N}$. These constraints induce an admissible family of candidate terms and, consequently, a prior-constrained dictionary
\[
\mathbf{\Theta}_{\mathcal{P}}=\big[\Phi^{(\mathcal{P})}_1(\mathbf{u}),\,\Phi^{(\mathcal{P})}_2(\mathbf{u}),\,\dots\big],
\]
where each basis function $\Phi_i^{(\mathcal{P})}$ implicitly satisfies the constraints in $\mathcal{P}$. Thus, while the five-stage workflow of Section~\ref{sec:pde_identification_framework} remains unchanged, Step~1 is replaced by a prior-constrained dictionary construction and Steps~2--5 operate on the resulting restricted dictionary. The approach applies to structural priors on $\mathcal{N}$. In Sections~\ref{subsec:hamiltonian_prior}–\ref{subsec:energy_prior}, we instantiate this design for three representative classes—Hamiltonian, conservation-law, and energy–dissipation priors—which exemplify geometric, flux-form (algebraic), and variational structures, respectively.

\noindent\textbf{Revised identification pipeline under prior constraints.}
\begin{itemize}
\item \textbf{Step $1^\prime$ (Prior-constrained dictionary construction).}
    Input: spatiotemporal measurements of the state $\mathbf{u}$ and structural constraints $\mathcal{P}$.
    Procedure: Construct the dictionary $\mathbf{\Theta}_{\mathcal{P}}$ by applying the structure-preserving operator (e.g., $J\nabla$, $\nabla \cdot$, or $\delta/\delta \mathbf{u}$) to a parameterized basis of the latent generating function (Hamiltonian, Flux, or Energy).
    Output: restricted dictionary $\mathbf{\Theta}_{\mathcal{P}}$ with tied parameters implied by $\mathcal{P}$.
    \item \textbf{Steps 2--5 (Sparse regression and model selection).}
    Apply the procedures of Section~\ref{sec:pde_identification_framework} to $(\mathbf{\Theta}_{\mathcal{P}},\,\partial_t \mathbf{u})$, using the same sparsity and selection criteria but acting on the restricted dictionary.
\end{itemize}

This physics-guided restriction provides three practical benefits.
(i) \emph{Improved conditioning.} Structural tying reduces the effective degrees of freedom by enforcing physical dependencies. This mitigates the ill-conditioning typical of generic overcomplete libraries by restricting the regression to a smaller subspace of admissible operators.
(ii) \emph{Noise robustness.} By excluding structurally inadmissible terms, the estimator’s variance is reduced and regularization is aligned with the true operator class, improving stability under measurement noise.
(iii) \emph{Interpretability.} Each retained term corresponds to a physically meaningful mechanism (e.g., a flux component), yielding models with transparent physical semantics.

We next demonstrate the prior-constrained dictionary construction using three classical structures: (a) \emph{Hamiltonian} systems (geometric/symplectic), (b) \emph{conservation-law} systems (algebraic flux-divergence form), and (c) \emph{energy–dissipation} systems (variational/gradient flows).
All three follow a common template: parameterize a latent physical object (Hamiltonian, flux, or energy), apply the structure-preserving operator, and generate PDE-level dictionary entries with coefficient tying induced by $\mathcal{P}$.

\subsection{Hamiltonian prior: skew-gradient structure}
\label{subsec:hamiltonian_prior}

For canonical variables $(q,p)\in\mathbb{R}^n\times\mathbb{R}^n$, Hamilton’s equations
\[
\dot{q}=\partial_p H(q,p),\qquad \dot{p}=-\partial_q H(q,p)
\]
can be written compactly by stacking $z=(q^\top,p^\top)^\top\in\mathbb{R}^{2n}$ as
\[
\dot{z}=J\nabla_z H(z),\qquad 
J=\begin{pmatrix}0&I_n\\-I_n&0\end{pmatrix},\quad J^\top=-J.
\]
The antisymmetry of $J$ implies energy conservation along trajectories: $\frac{d}{dt}H(z(t))=\nabla H\cdot \dot z=\nabla H\cdot J\nabla H=0$. Moreover, Liouville’s theorem gives a divergence-free phase flow since $\nabla\!\cdot(J\nabla H)=\mathrm{tr}(J\nabla^2 H)=0$ for constant $J$. These identities encode the underlying symplectic geometry. (Field-theoretic Hamiltonian PDEs can be treated analogously by replacing $J$ with a suitable Poisson operator.)

\paragraph{Dictionary design}
Let $\{\phi_\alpha(q,p)\}_{\alpha\in\mathcal{A}}$ be a basis for $H$ and parameterize
\[
H(q,p)=\sum_{\alpha\in\mathcal{A}} w_\alpha\,\phi_\alpha(q,p),\qquad w_\alpha\in\mathbb{R}.
\]
Differentiation yields the PDE-level dictionary as skew-gradient dictionary entries
\[
\mathcal{D}_{\mathrm{Ham}}
=\Big\{\, J\nabla \phi_\alpha(q,p)\ \Big|\ \alpha\in\mathcal{A}\Big\},\qquad
\dot z=\sum_{\alpha\in\mathcal{A}} w_\alpha\, J\nabla \phi_\alpha(q,p).
\]
The coefficients $\{w_\alpha\}$ are \emph{shared} across state components: the same $w_\alpha$ multiplies both the $q$- and $p$-components through $J\nabla\phi_\alpha$, guaranteeing the skew-gradient structure by construction.

This parameter tying reduces the number of free parameters and enforces energy conservation and volume preservation in the identified model.

\subsection{Conservation-law prior: flux-form representation}
\label{subsec:conservation_prior}
Many continuum models admit a conservative form
\[
\mathbf{u}_t+\nabla\!\cdot F(\mathbf{u})=0,\qquad \mathbf{u}:\Omega\to\mathbb{R}^n,
\]
which states that the temporal change of $\mathbf{u}$ in a control volume is balanced by the net outward flux. Under periodic or no-flux boundary conditions, integrating over $\Omega$ yields conservation of the corresponding integral invariants. The scalar case $n=1$ is recovered by writing $u_t+\nabla\!\cdot F(u)=0$.

\paragraph{Dictionary design}
Let $\{\phi_\alpha(\mathbf{u})\}_{\alpha\in\mathcal{A}}$ be a basis for the flux $F(\mathbf{u})$ and parameterize
\[
F(\mathbf{u})=\sum_{\alpha\in\mathcal{A}} w_\alpha\,\phi_\alpha(\mathbf{u}),\qquad w_\alpha\in\mathbb{R}.
\]
Differentiation yields the PDE-level dictionary as divergence dictionary entries
\[
\mathcal{D}_{\mathrm{Flux}}
=\Big\{\, \nabla\cdot \phi_\alpha(\mathbf{u})\ \Big|\ \alpha\in\mathcal{A}\Big\},\qquad
\partial_t \mathbf{u} + \sum_{\alpha\in\mathcal{A}} w_\alpha\, \nabla\cdot \phi_\alpha(\mathbf{u}) = 0.
\]

Coefficient tying encodes physical symmetries; e.g., for isotropic terms in 2D, enforce $w_{x,\alpha}=w_{y,\alpha}$ across directional flux components. In one dimension we keep the conservative notation to emphasize structure-preserving discretization,
\[
\mathcal{D}_{\mathrm{Flux}}^{1\mathrm{D}}
=\{(u)_x,\ (u^2)_x,\ (u_x)_x,\ (u^2 u_x)_x,\ \dots\},
\]
which avoids confusion with non-conservative derivative libraries while making the underlying flux explicit. For numerical implementation, the continuous divergence should be paired with its discrete counterpart (finite-volume or compatible finite-difference operators) to preserve conservation at the discrete level.

\subsection{Energy-dissipation prior: gradient-flow structure}
\label{subsec:energy_prior}
A broad class of dissipative dynamics can be written as gradient flows driven by an energy functional,
\[
\mathbf{u}_t = -\,\frac{\delta E}{\delta \mathbf{u}}, 
\qquad 
E[\mathbf{u}]=\int_\Omega \mathcal{E}\!\big(\mathbf{u},\nabla \mathbf{u},\nabla^2 \mathbf{u},\dots\big)\,dx,
\]
with $\mathbf{u}:\Omega\to\mathbb{R}^m$. Under periodic or no-flux boundary conditions,
\[
\frac{dE}{dt}
=\Big\langle \frac{\delta E}{\delta \mathbf{u}},\,\mathbf{u}_t\Big\rangle
= -\Big\|\frac{\delta E}{\delta \mathbf{u}}\Big\|_{L^2(\Omega)}^2 \le 0,
\]
so $E$ decays monotonically.
Canonical examples include the heat equation ($\mathcal{E}= \tfrac{1}{2}\lvert\nabla u\rvert^2$ yielding $u_t=\Delta u$), Allen--Cahn ($\mathcal{E}=\tfrac{\epsilon^2}{2}\lvert\nabla u\rvert^2+\tfrac{1}{4}(u^2-1)^2$ yielding $u_t=\epsilon^2\Delta u - (u^3-u)$)

We adopt this representation as a prior on the operator class, thereby restricting candidate terms to variational derivatives of admissible energy densities.
The vector-valued case is handled natively by the formulation above: $\mathbf{u}\in\mathbb{R}^n$ and $\delta E/\delta \mathbf{u}$ is taken componentwise. We present the general vector framework, which encompasses the scalar setting as the special case $n=1$. 

\paragraph{Dictionary design}
Let $\{\phi_\alpha(\mathbf{u},\nabla \mathbf{u},\dots)\}_{\alpha\in\mathcal{A}}$ be a basis for the admissible energy density $\mathcal{E}$ and parameterize the total energy functional as
\[
E[\mathbf{u}] = \int_\Omega \sum_{\alpha\in\mathcal{A}} w_\alpha\,\phi_\alpha(\mathbf{u},\nabla \mathbf{u},\dots)\, dx, \qquad w_\alpha\in\mathbb{R}.
\]
Variational differentiation yields the PDE-level dictionary as gradient-flow entries
\[
\mathcal{D}_{\mathrm{GF}}
=\Big\{\, -\frac{\delta}{\delta \mathbf{u}}\int_\Omega \phi_\alpha\,dx \ \Big|\ \alpha\in\mathcal{A}\Big\},
\qquad
\partial_t \mathbf{u} = \sum_{\alpha\in\mathcal{A}} w_\alpha \left( -\frac{\delta}{\delta \mathbf{u}} \int_\Omega \phi_\alpha \, dx \right).
\]

For numerical fidelity, the discrete variational derivative should be consistent with the discrete energy (e.g., finite-element/finite-difference pairs that satisfy a discrete chain rule), ensuring that the learned model preserves the intended dissipation at the discrete level.

\section{Case studies}
\label{sec:case_studies}

In this section, we demonstrate the proposed prior-informed dictionary refinement through case studies across the three structural classes discussed in the previous section.

\noindent\textbf{Weak formulation setup.} \label{sec:case_studies:weak}
Note that all three cases mentioned above admit weak formulations. To mitigate measurement noise, we adopt a unified weak formulation. Consider a general system $u_t - \mathcal{N}(u, \dots)=0$. Let $\psi\in C_c^\infty(\Omega\times[0,T))$ be a smooth test function with compact temporal support and appropriate spatial boundary conditions (periodic or homogeneous). Multiplying the system by $\psi$, integrating over $\Omega\times[0,T]$, and integrating by parts yield the weak identity
\begin{equation}
\int_0^T\!\!\int_\Omega u\,\mathcal{F}^*[\psi]\,dx\,dt=0,
\end{equation}
where $\mathcal{F}^*$ is the formal adjoint operator. This transfers derivatives from the noisy state $u$ to the smooth test function $\psi$.

\subsection{Hamiltonian prior: skew-gradient structure}
\label{subsec:case_hamiltonian}

We begin with Hamiltonian systems, characterized by energy conservation and symplectic phase-space dynamics. We examine two canonical examples: the harmonic oscillator and the three-body problem.

\noindent\textbf{(a) Harmonic oscillator.}

The harmonic oscillator serves as a canonical model of coupled position–momentum dynamics. For a specific parameter configuration, the governing equations are:\begin{equation}\dot{q} = 2p, \qquad \dot{p} = -2q,\label{eq:HO_nonham}\end{equation}where $q(t)$ and $p(t)$ denote displacement and momentum, respectively.

Consistent with Section~\ref{subsec:hamiltonian_prior}, the system admits a Hamiltonian formulation with state $z = (q, p)^\top$. The dynamics are governed by the total energy:\begin{equation}H(q,p) = \frac{p^2}{2m} + \frac{1}{2}k q^2 = p^2+q^2,\label{eq:HO_Hamiltonian}\end{equation}where $m=0.5$ represents the mass and $k=2$ denotes the spring stiffness. The system satisfies the canonical evolution $\dot{z} = J \nabla_z H(z)$.

\noindent\textbf{Symplectic constraint for identification.}
The existence of the scalar potential \eqref{eq:HO_Hamiltonian} implies that the vector field is strictly determined by the symplectic gradient structure. In Section \ref{numericalexp:harmonic_spr}, we utilize this prior to constrain the identification to energy-conserving terms.

\noindent\textbf{Weak formulation under the Hamiltonian prior.}
Let \(\boldsymbol{\psi}(q,p,t)\in C_c^\infty(\Omega\times(0,T);\mathbb{R}^2)\) be a smooth, compactly supported test function. 
Multiplying \(\dot{z}=J\nabla_z H(z)\) by \(\boldsymbol{\psi}\) and integrating over \(\Omega\times[0,T]\) yield
\begin{equation}
-\int_0^T\!\!\int_{\Omega} z\cdot \boldsymbol{\psi}_t \,dq\,dp\,dt
-\int_0^T\!\!\int_{\Omega} H(q,p)\,\nabla_z\!\cdot\!\big(J\boldsymbol{\psi}\big)\,dq\,dp\,dt = 0.
\end{equation}

\noindent\textbf{(b) Three-Body problem.}  
The three-body problem describes three particles interacting under Newtonian gravitation in three-dimensional space. Let $\mathbf{q}_i, \mathbf{p}_i \in \mathbb{R}^3$ denote the position and momentum of the $i$-th body.
The equations of motion are
\begin{equation}
\begin{aligned}
\dot{\mathbf{q}}_i &= \mathbf{p}_i / m_i, \\[2pt]
\dot{\mathbf{p}}_i &= -\sum_{j\neq i} G m_i m_j \frac{\mathbf{q}_i - \mathbf{q}_j}{\|\mathbf{q}_i - \mathbf{q}_j\|^3},
\quad i=1,2,3,
\end{aligned}
\label{eq:threebody_nonham}
\end{equation}
where $G$ is the gravitational constant and $m_i$ is the mass. This system exhibits highly coupled, nonlinear, and potentially chaotic dynamics.

Defining the canonical state $\mathbf{z} = (\mathbf{q}_1^\top, \dots, \mathbf{q}_3^\top, \mathbf{p}_1^\top, \dots, \mathbf{p}_3^\top)^\top \in \mathbb{R}^{18}$, the system admits the Hamiltonian form $\dot{\mathbf{z}} = J \nabla_{\mathbf{z}} H(\mathbf{z})$, governed by the total energy:
\begin{equation}
H(\mathbf{q},\mathbf{p})= \sum_{i=1}^{3} \frac{\|\mathbf{p}_i\|^2}{2m_i} - \sum_{i<j} \frac{Gm_i m_j}{\|\mathbf{q}_i - \mathbf{q}_j\|}.
\label{eq:threebody_H}
\end{equation}

This formulation ensures conservation of energy and angular momentum. Under the \emph{Hamiltonian prior}, we identify the scalar Hamiltonian $H(\mathbf{q},\mathbf{p})$ rather than the component-wise equations. 
We parameterize $H(\mathbf{q},\mathbf{p})$ using weights $w_\alpha$:
\begin{equation} \label{eq:threebody_basis}
H(\mathbf{q},\mathbf{p}) = \sum_{\alpha=1}^{M} w_\alpha h_\alpha (\mathbf{q},\mathbf{p}),
\end{equation}
yielding the dictionary of symplectic gradients $\mathcal{D}_{\text{ham}} = \{ J\nabla h_\alpha \}_{\alpha=1}^{M}$.

\noindent\textbf{Symplectic constraint for identification.}
This formulation ensures conservation of energy and angular momentum. In Section \ref{numericalexp:3bdy}, we leverage this Hamiltonian structure to reduce the learning problem from identifying 18 coupled equations to identifying a single scalar energy functional $H(\mathbf{z})$.

\noindent\textbf{Weak formulation under the Hamiltonian prior.}
Let $\boldsymbol{\psi}(\mathbf{z},t)\in C_c^\infty(\mathcal{Z}\times(0,T);\mathbb{R}^{18})$ be a test function on the $18$-dimensional phase space $\mathcal{Z}$.

The weak formulation yields
\begin{equation}
-\int_0^T\!\!\int_{\mathcal{Z}}
\mathbf{z}\cdot\boldsymbol{\psi}_t\,d\mathbf{z}\,dt
-\int_0^T\!\!\int_{\mathcal{Z}}
H(\mathbf{z})\,\nabla_{\mathbf{z}}\!\cdot\!\big(J\boldsymbol{\psi}\big)\,d\mathbf{z}\,dt
=0.
\end{equation}

\subsection{Conservation-Law prior: flux form representation}
\label{subsec:case_conservation}

We consider conservation laws, where temporal variation is balanced by flux divergence. Two benchmark PDEs are examined: the 1D Burgers' equation and the 2D shallow-water equations.

\noindent\textbf{(a) 1D Burgers' equation.}  
The one-dimensional Burgers' equation models nonlinear advective transport of a scalar field $u(x,t)$ and admits the equivalent forms
\begin{equation}
u_t + u\,u_x = 0, 
\qquad
u_t + \bigl(\tfrac{1}{2}u^2\bigr)_x = 0,
\label{eq:Burgers_conservative}
\end{equation}
where the conservative form highlights the flux balance with conservative variable $q=u$ and flux $F(q)=\tfrac{1}{2}q^2$, i.e.,
\begin{equation}
q_t + F(q)_x = 0.
\label{eq:Burgers_flux_form}
\end{equation}
This representation supports the \emph{conservation-law prior}, restricting admissible features to spatial derivatives of fluxes.

\noindent\textbf{Conservation constraint for identification.}
The form \eqref{eq:Burgers_flux_form} implies that the dynamics are strictly governed by the divergence of a flux. In Section \ref{numerical_exp:burgers}, we enforce this prior by restricting the candidate library to terms of the form $(F(u))_x$.

\noindent\textbf{Weak formulation under the conservation-law prior.}

Consider a test function $\psi(x,t) \in C_c^\infty(\Omega \times [0,T))$. By multiplying \eqref{eq:Burgers_flux_form} by $\psi$ and integrating by parts over the spatiotemporal domain, the weak identity is given by
\begin{equation}
-\int_0^T\!\!\int_\Omega q\,\psi_t\,dx\,dt - \int_0^T\!\!\int_\Omega F(q)\,\psi_x\,dx\,dt = 0,
\label{eq:Burgers_weak_general}
\end{equation}
where $F(q)$ is a sparse linear combination of terms from the candidate flux library.

\noindent\textbf{(b) 2D Shallow-Water equations.}
The two-dimensional shallow-water (SW) equations describe an incompressible fluid layer of depth $h(x,y,t)$ with horizontal velocity $(u,v)$. The governing equations are
\begin{equation}
\begin{aligned}
h_t + (hu)_x + (hv)_y &= 0,\\
u_t + u\,u_x + v\,u_y + g\,h_x &= 0,\\
v_t + u\,v_x + v\,v_y + g\,h_y &= 0,
\end{aligned}
\label{eq:SWE_nonconservative}
\end{equation}
where $g$ denotes gravitational acceleration. Here, the first equation enforces mass continuity, while the latter two describe the conservation of horizontal momentum.

Define the conservative state $\mathbf{q}=(h,hu,hv)^\top$. The system admits a vector form of conservation-law:
\begin{equation}
\mathbf{q}_t + \mathbf{F}(\mathbf{q})_x + \mathbf{G}(\mathbf{q})_y = 0,
\label{eq:swe_flux}
\end{equation}
where $\mathbf{F}(\mathbf{q})$ and $\mathbf{G}(\mathbf{q})$ represent the unknown nonlinear flux functions in the $x$- and $y$-directions, respectively. In our framework, these fluxes are treated as latent physical objects to be identified from an overcomplete dictionary of candidate functional terms.

\noindent\textbf{Conservation constraint for identification.} The structural identity \eqref{eq:swe_flux} implies that the divergence of admissible fluxes strictly governs the dynamics. Following the prior-informed strategy, we do not regress the components of $\mathbf{q}_t$ independently. Instead, we construct a physically consistent dictionary where each candidate feature is a divergence term $\nabla \cdot (\mathbf{F}, \mathbf{G})^{\top}$, ensuring that any identified model preserves the fundamental mass and momentum conservation laws by construction.

\noindent\textbf{Weak formulation under the conservation-law prior.} Consider a smooth test function $\boldsymbol{\psi}(x,y,t)\in C_c^\infty(\Omega\times[0,T);\mathbb{R}^3)$. By transferring the spatial derivatives onto the test function, the weak formulation of \eqref{eq:swe_flux} is given by:\begin{equation}-\int_0^T\!\!\int_{\Omega} \mathbf{q} \cdot \boldsymbol{\psi}_t \,dx\,dy\,dt - \int_0^T\!\!\int_{\Omega} \big(\mathbf{F}(\mathbf{q}) \cdot \boldsymbol{\psi}_x + \mathbf{G}(\mathbf{q}) \cdot \boldsymbol{\psi}_y \big) \,dx\,dy\,dt = 0.\label{eq:SWE_weak_general}\end{equation}

\subsection{Energy-Dissipation prior: gradient flow structure}
\label{subsec:case_energy}

Finally, we consider systems evolving to minimize an energy functional, exhibiting dissipative behavior. Two canonical examples are the 1D diffusion equation and the 1D Allen--Cahn equation. 

\noindent\textbf{(a) 1D Diffusion equation.}\label{subsec:case_energy_diffusioneq}
The one-dimensional diffusion equation describes the dissipative spreading of a scalar field $u(x,t)$ on $\Omega=(a,b)$, serving as a canonical prototype for gradient flows.
Its standard form is
\begin{equation}
u_t = \nu\,u_{xx},
\qquad \nu>0,
\label{eq:diffusion_strong}
\end{equation}
where $\nu$ is the diffusivity coefficient.

This equation is the $L^2$-gradient flow of the \emph{Dirichlet energy functional}:
\begin{equation}
E[u] = \frac{\nu}{2} \int_{\Omega} |u_x|^2 \,dx.
\label{eq:diffusion_energy_def}
\end{equation}

The associated gradient-flow dynamics recover the diffusion equation~\eqref{eq:diffusion_strong}:
\begin{equation}
u_t = -\frac{\delta E}{\delta u} = \nu\,u_{xx}.
\label{eq:diffusion_gradient_flow}
\end{equation}

This formulation highlights the intrinsic \emph{energy-dissipation structure}, where temporal evolution follows the steepest descent of $E[u]$.

\noindent\textbf{Gradient flow constraint for identification.}
Since the diffusion equation is formally the $L^2$-gradient flow of the Dirichlet energy \eqref{eq:diffusion_energy_def}, any data-driven identification should ideally respect this variational structure. In our experiments (Section \ref{numericalexp:diffu}), we will leverage this property by parameterizing the energy density rather than the PDE terms directly.

\noindent\textbf{Weak formulation under the energy–dissipation prior.}

Consider a smooth test function $\psi(x,t) \in C_c^\infty(\Omega \times [0,T))$. Multiplying the gradient-flow identity $u_t = -\frac{\delta E}{\delta u}$ by $\psi$ and integrating by parts over the spatiotemporal domain yield:
\begin{equation}
-\int_0^T\!\!\int_\Omega u\,\psi_t\,dx\,dt + \int_0^T\!\!\int_\Omega \frac{\delta E}{\delta u} \psi \,dx\,dt = 0,
\label{eq:diffusion_weak_general}
\end{equation}
where $\frac{\delta E}{\delta u}$ is the variational derivative of the latent energy functional $E[u]$.

\noindent\textbf{(b) 1D Allen--Cahn equation.} \label{subsec:case_energy_allencahn}
The one-dimensional Allen--Cahn equation models the relaxation dynamics of a scalar order parameter $u(x,t)$ on $\Omega=(a,b)$, serving as a canonical prototype for phase separation and interface motion.
Its standard form is
\begin{equation}
u_t = u_{xx} - (u^3 - u),
\label{eq:allen_cahn_strong}
\end{equation}
where $u_{xx}$ enforces smoothness and $-(u^3 - u)$ drives the system toward stable equilibria $u = \pm 1$.

This equation is the $L^2$-gradient flow of the \emph{double-well energy functional}:
\begin{equation}
E[u] = \int_{\Omega}
\Big( \tfrac{1}{2}|u_x|^2 + \tfrac{1}{4}(u^2 - 1)^2 \Big) \,dx.
\label{eq:allen_cahn_energy_def}
\end{equation}
where the first term penalizes spatial variations and the second favors the symmetric minima $u=\pm1$.

The associated gradient flow recovers \eqref{eq:allen_cahn_strong}:
\begin{equation}
u_t
= -\frac{\delta E}{\delta u}
= u_{xx} - (u^3 - u).
\label{eq:allen_cahn_gradient_flow}
\end{equation}

\noindent\textbf{Gradient flow constraint for identification.}
The variational structure above implies that the admissible driving forces in the PDE must originate from a bounded energy potential. In Section \ref{numericalexp:allen_cahn}, we leverage this property to restrict the candidate library, ensuring that the identified model satisfies the energy dissipation law.

\noindent\textbf{Weak formulation under the energy–dissipation prior (Allen–Cahn).}

Consider a smooth test function $\psi(x,t) \in C_c^\infty(\Omega \times [0,T))$. Following the variational structure established in \eqref{eq:allen_cahn_gradient_flow}, the weak formulation is expressed as:
\begin{equation}
-\int_0^T\!\!\int_\Omega u\,\psi_t\,dx\,dt + \int_0^T\!\!\int_\Omega \frac{\delta E}{\delta u} \psi \,dx\,dt = 0,
\label{eq:AC_weak_general}
\end{equation}
where $\frac{\delta E}{\delta u}$ is the variational derivative of the unknown energy functional $E[u]$.

\section{Numerical experiment} \label{sec:numericalexp}
With the prior informed dictionary constructions, this section evaluates identification quality and robustness to measurement noise for several PDEs. For each example, we compare the following four configurations:
\begin{itemize}
    \item[Conf. 1] Separate identification without weak form: Directly identify PDEs without using prior nor weak formulation, corresponding to the method in \cite{he2022robust} without Successively Denoised Differentiation (SDD).
    \item[Conf. 2] Separate identification with weak form: Directly identify PDEs without using prior but use weak formulation, corresponding to the method in \cite{tang2022weakident,TANG2026114454}.
    \item[Conf. 3] Prior-based identification without weak form: Use prior to identify a PDE without using weak formulation.
    \item[Conf. 4] Prior-based identification with weak form: The proposed method.
\end{itemize}
For each method, we measure its performance by the TPR defined as: 
    \begin{equation} \label{tpr}
        \text{TPR} = \frac{|\{ l : \textbf{c}^*(l) \neq 0 \ \text{and} \ \textbf{c}(l) \neq 0 \}|}{|\{ l : \textbf{c}^*(l) \neq 0 \}|}.
    \end{equation}
    This metric computes the ratio of correctly identified non-zero coefficients to the total number of true non-zero coefficients. It quantifies the method's ability to correctly identify relevant features in the underlying differential equation.

Consider a space-time domain $\Omega\times [0,T]$ where $\Omega = [L_1,L_2]$ with $0\leq L_1<L_2$ and $T>0$. The given data (discretized PDE solution with noise) is denoted by
\begin{align}
\mathcal{D}=\{U_i^n = u(x_i,t^n)+\varepsilon_i^n, i = 0, 1,\ldots,N^x, n = 0, 1,\ldots,N^t\}\label{dataset},
\end{align}
where $x_i=i\Delta x \in [L_1,L_2] $ and $t^n=n\Delta t\in [0,T]$; $\varepsilon_i^n$ represents independent zero-mean noise; and $u(x_i,t^n)$ represents the clean data, for $i = 0, 1,\ldots,N^x,$ and $ n = 0, 1,\ldots,N^t$.

Let \[U_{\text{max}} = \max_{j=1,\ldots,N^x} \max_{m=1,\ldots,N^t} U_j^m, \quad U_{\text{min}} = \min_{j=1,\ldots,N^x} \min_{m=1,\ldots,N^t} U_j^m.\] We consider a noise model $\epsilon\sim\mathcal{N}(0,\sigma^2)$ with  standard deviation \(\sigma\)  defined as:
\[
\sigma = \frac{\sigma_{\text{NSR}}}{N^t N^x} \sum_{j=1}^{N^x} \sum_{m=1}^{N^t} \left| U_j^m - \frac{U_{\text{max}} + U_{\text{min}}}{2} \right|^2.
\]
Here, $\sigma_{\text{NSR}}\in [0,1]$ indicates the noise-to-signal ratio (NSR).

We also visually compare the true and identified trajectories across different noise levels, as illustrated in Figures~\ref{fig:threebody_identified_trajectory_all_3D},~\ref{fig:SWE_h_trajectory_grid_00new}--\ref{fig:SWE_h_trajectory_grid_50new},~\ref{fig:AllenCahn_2x3_grid}. Here, `true trajectories' refer to the ground truth data, while `identified trajectories' represent the solutions obtained from the identified PDEs/ODEs starting from the same initial conditions.

\subsection{Hamiltonian prior}
\subsubsection{Harmonic oscillator} \label{numericalexp:harmonic_spr}
We begin with the harmonic oscillator, serving as a benchmark for Hamiltonian identification (see Section \ref{subsec:case_hamiltonian}). The data is generated on a uniform grid covering the phase-space domain $[-1.2, 1.2]^2$. Trajectories are integrated over the interval $t \in [0, 3]$ with a time step $\Delta t = 0.01$, initialized from energy level sets with radii $r \in [0.1, 1.0]$. Given the Hamiltonian $H = p^2 + q^2$, these radii correspond to concentric circular orbits in phase space with total energy $E = r^2$.

\medskip
\noindent\textbf{Feature libraries.}
We construct two feature libraries to evaluate the Hamiltonian prior: a generic polynomial library and a Hamiltonian-constrained library.

\noindent\textbf{Baseline (non-prior) dictionary.}
We construct a baseline polynomial library $\mathcal{D}_{\text{base}}$ using a basis $\mathcal{B}$ up to third order:$$\mathcal{B} = \{ 1, p, q, pq, p^2, q^2, p^2q, pq^2,p^3,q^3 \}.$$

\noindent\textbf{Hamiltonian prior dictionary.}
Based on the symplectic structure discussed in Section \ref{subsec:case_hamiltonian}, the prior-informed dictionary is generated by applying the skew-gradient operator to the polynomial basis:
\begin{equation}\label{HO_dic_ham}
\mathcal{D}_{\text{Ham}} = \big\{ J\nabla_z \phi \mid \phi \in \mathcal{B} \big\}.
\end{equation}
This ensures every candidate model in $\mathcal{D}_{\text{Ham}}$ inherently satisfies energy conservation.

\begin{figure}[ht!]
\centering
\begin{subfigure}{0.45\linewidth}
  \includegraphics[width=\linewidth]{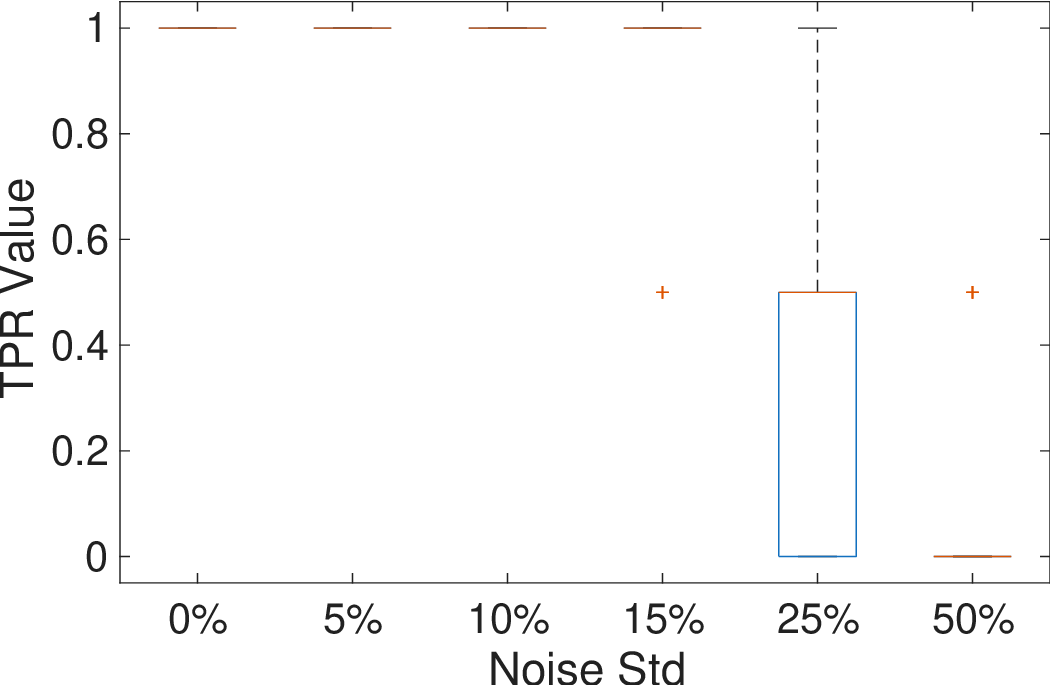}
  \caption{Separate (Strong form)} 
  \label{subfig:sp_separate_spring}
\end{subfigure} 
\hfill
\begin{subfigure}{0.45\linewidth}
  \includegraphics[width=\linewidth]{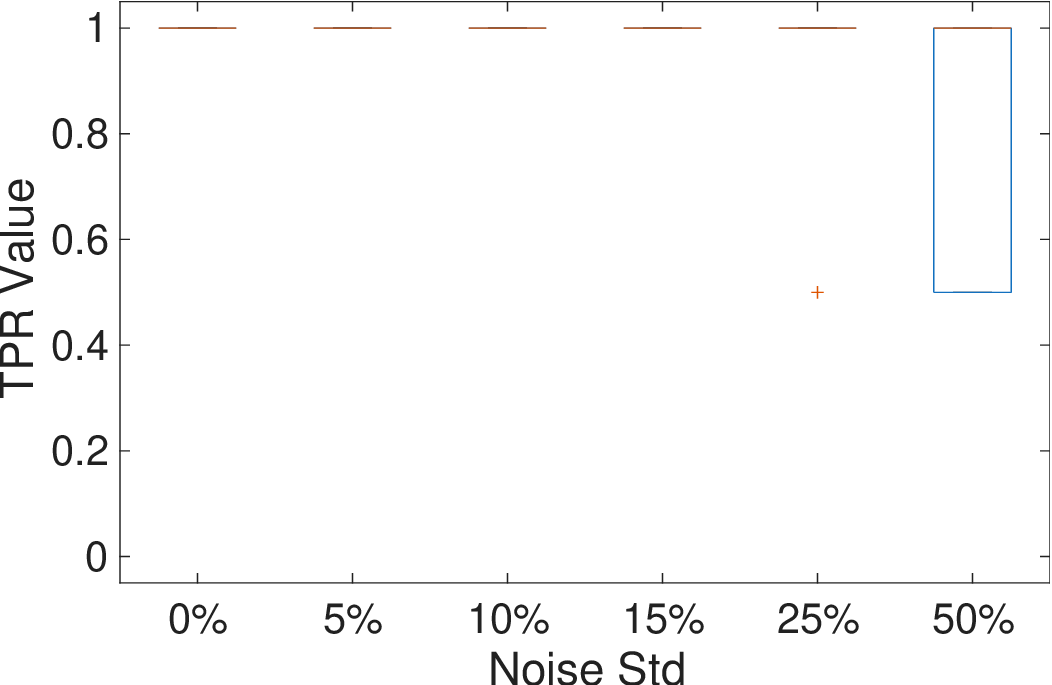}
  \caption{Separate (Weak form)}  
  \label{subfig:weak_separate}
\end{subfigure}

\vspace{-0em}

\begin{subfigure}{0.45\linewidth}
  \includegraphics[width=\linewidth]{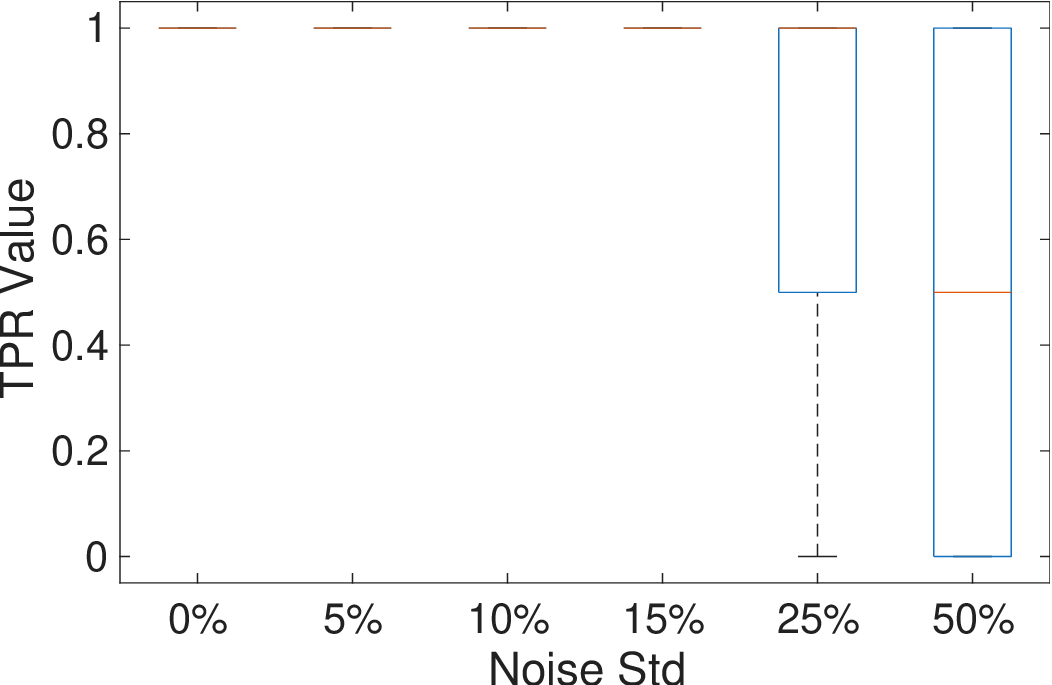}
  \caption{Prior-based (Strong form)}  
  \label{subfig:sp_prior_spring}
\end{subfigure}
\hfill
\begin{subfigure}{0.45\linewidth}
  \includegraphics[width=\linewidth]{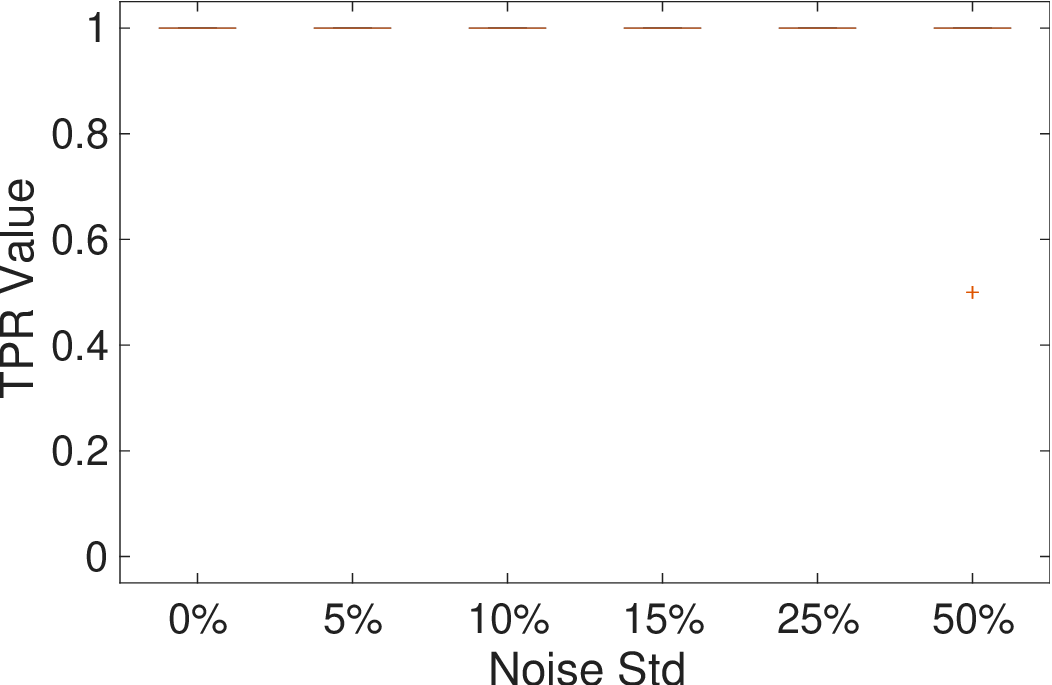}
  \caption{Prior-based (Weak form)} 
  \label{subfig:weak_prior}
\end{subfigure}

\vspace{-0em}

\caption{
TPR results from 20 repeated identification trials for the harmonic oscillator system with different noise levels $\{0\%, 5\%, 10\%, 15\%, 25\%, 50\% \}$ under four configurations as mentioned above.
}
\label{fig:Hamiltonian_system}
\end{figure}

\noindent\textbf{Identification results and performance.}
We identify the PDE system using both the baseline and prior-informed libraries with the four configurations. Under various noise levels, the TPR averaged over 20 trials are shown in Figure \ref{fig:Hamiltonian_system}. With the baseline library, the position and momentum equations are identified independently. In contrast, the Hamiltonian prior constrains the search to a single joint regression, reducing the degree of freedom. Across all noise levels, the proposed method has stable performances. It correctly identifies the active support $\mathcal{S}^*_{\text{Ham}}$ corresponding to the quadratic energy components $p^2$ and $q^2$ with TPR=1, as shown in (d). While results with other configurations have obvious reduced TPR with high noise levels. The results demonstrate improved performance with our proposed method: the Hamiltonian prior ensures energy-preserving, physically consistent recovery, while the weak-form integration substantially enhances numerical stability and robustness to noise.

\subsubsection{Three-Body problem} \label{numericalexp:3bdy}

We test the proposed method on the three-body problem described in Section \ref{subsec:case_hamiltonian}. The system dynamics are governed by the Hamiltonian $H(\mathbf{z})$ defined in Equation \eqref{eq:threebody_H}. In our experiments, we adopt normalized units with gravitational constant $G=1.0$ and equal masses $m_i=1.0$. The system is initialized in a perturbed figure-eight configuration, with positions $\mathbf{q}_1 = -\mathbf{q}_3 = (-0.9700, 0.2431, 0)^\top$ and $\mathbf{q}_2 = \mathbf{0}$. Initial momenta are set to $\mathbf{p}_1 = (0.4662, 0.4324,$\\$ 0.001)^\top$, $\mathbf{p}_2 = (-0.9324, -0.8647, 0)^\top$, and $\mathbf{p}_3 = (0.4662, 0.4324, -0.001)^\top$, ensuring a stationary center of mass and quasi-planar motion. Trajectories are integrated over $t \in [0, 100]$ with a time step $\Delta t = 0.01$, confined within the domain $[-1.2, 1.2]^3$.

\medskip
\noindent\textbf{Feature libraries.}
We construct two feature libraries and compare an unconstrained dictionary with a Hamiltonian-constrained dictionary enforcing symplectic structure.

\noindent\textbf{Baseline (non-prior) dictionary.}
The baseline library combines single-body polynomials and pairwise gravitational terms to regress $\dot{\mathbf{q}}_i$ and $\dot{\mathbf{p}}_i$:
\begin{equation}
\mathcal{D}_{\text{base}} = \Big\{ \mathbf{q}_i, \mathbf{p}_i, \mathbf{q}_i^2, \mathbf{p}_i^2, \mathbf{q}_i \mathbf{p}_i, \|\mathbf{q}_{ij}\|^{-1}, \mathbf{q}_{ij}\|\mathbf{q}_{ij}\|^{-3} \Big\}_{i,j \in \{1,2,3\}, j \neq i},
\label{eq:3body_baseline_dic}
\end{equation}
where $\mathbf{q}_{ij} = \mathbf{q}_i - \mathbf{q}_j$ denotes the relative position vector between bodies $i$ and $j$. This library treats each component of the canonical state $\mathbf{z}=(\mathbf{q}, \mathbf{p})^\top$ independently, thereby ignoring the underlying physical symmetries, conservation laws, and symplectic coupling.

\noindent\textbf{Hamiltonian prior dictionary.}
We restrict admissible terms to symplectic gradients of a scalar energy. To validate the recovery of isotropy from data, we construct the energy basis $\mathcal{B}_{H}$ using scalar momentum components and pairwise potentials:
\begin{equation}
\mathcal{B}_{H} = \Big\{ p_{i,\alpha}^2, \|\mathbf{q}_{ij}\|^{-1} \mid i,j \in \{1,2,3\}, j > i, \alpha \in \{x,y,z\} \Big\},
\label{eq:3body_ham_basis}
\end{equation}
where $p_{i,\alpha}$ denotes the $\alpha$-th component of the momentum for body $i$. Applying the symplectic gradient to each basis element yields the Hamiltonian dictionary:
\begin{equation}
\mathcal{D}_{\text{Ham}} = \big\{ J \nabla_{\mathbf{z}}\phi \mid \phi \in \mathcal{B}_{H} \big\}.
\label{eq:3body_ham_dic}
\end{equation}
This construction ensures that every candidate model naturally respects the symplectic structure of the system.

\begin{figure}[ht!]
\centering
\begin{subfigure}{0.45\linewidth}
  \includegraphics[width=\linewidth]{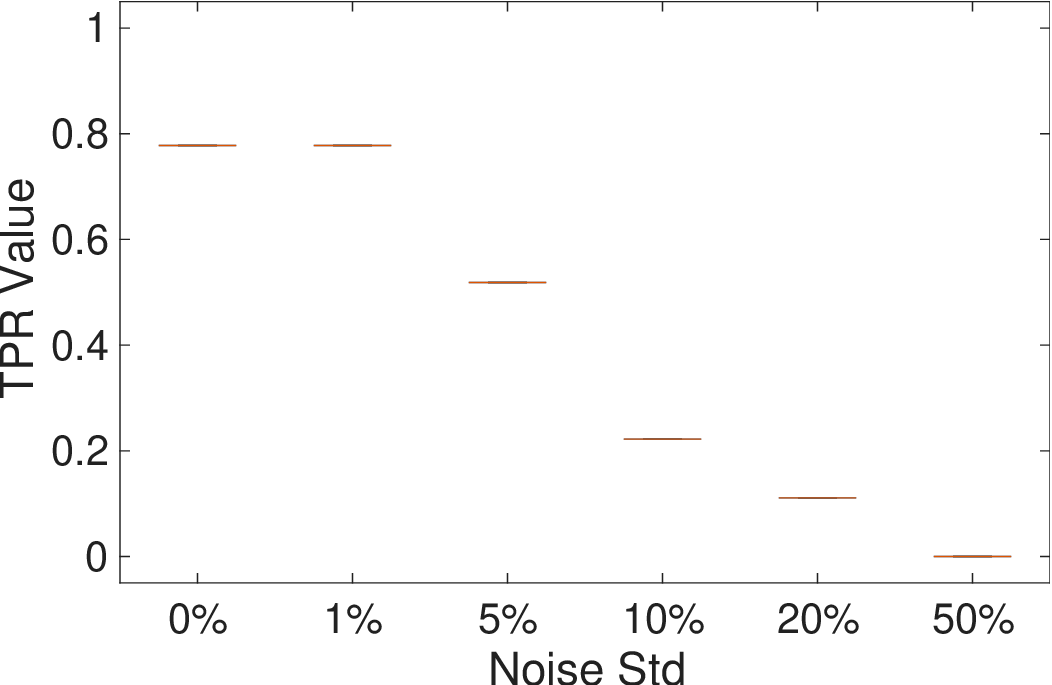}
  \caption{Separate (Strong form)}
  \label{subfig:prior_sp}
\end{subfigure}
\hfill
\begin{subfigure}{0.45\linewidth}
  \includegraphics[width=\linewidth]{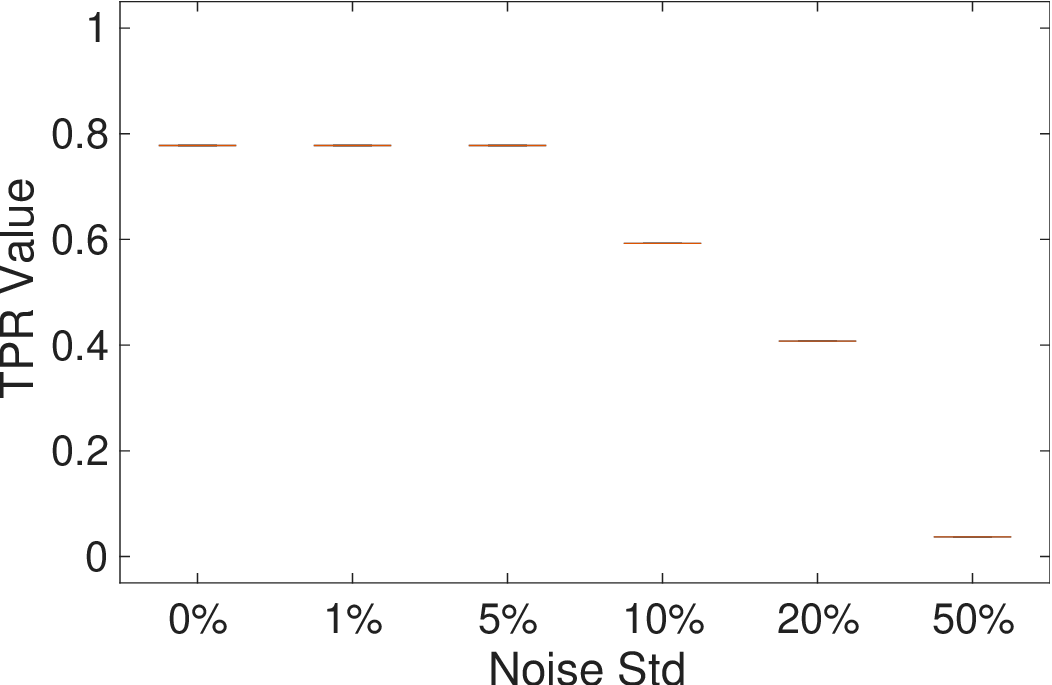}
  \caption{Separate (Weak form)}
  \label{subfig:prior_weak}
\end{subfigure}

\vspace{-0em}

\begin{subfigure}{0.45\linewidth}
  \includegraphics[width=\linewidth]{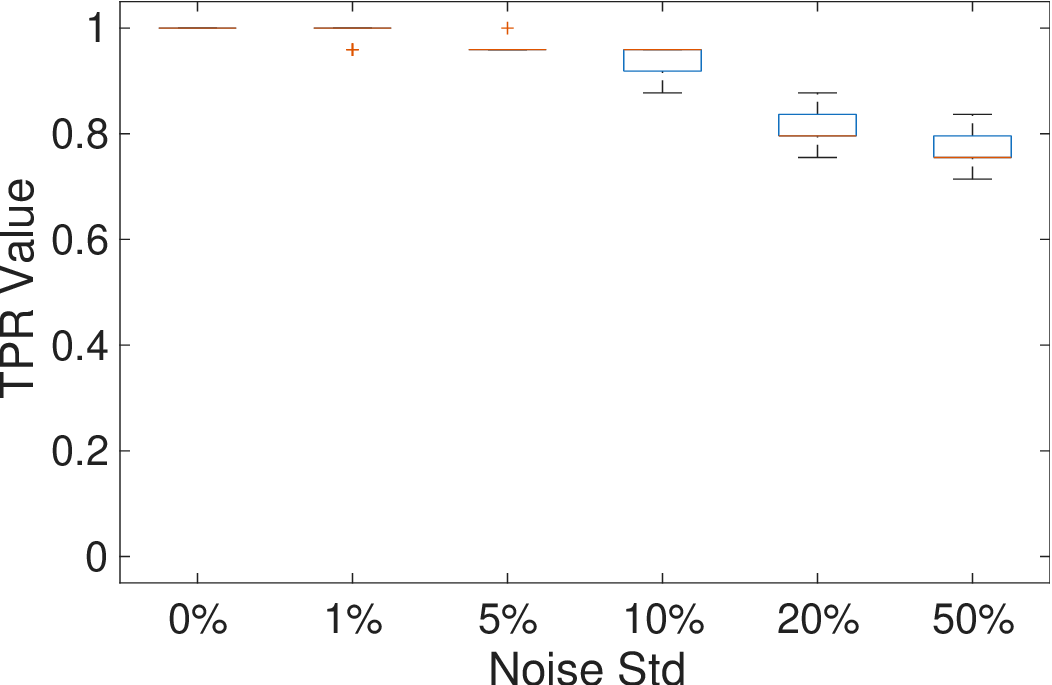}
  \caption{Prior-based (Strong form)}
  \label{subfig:no_prior_sp}
\end{subfigure}
\hfill
\begin{subfigure}{0.45\linewidth}
  \includegraphics[width=\linewidth]{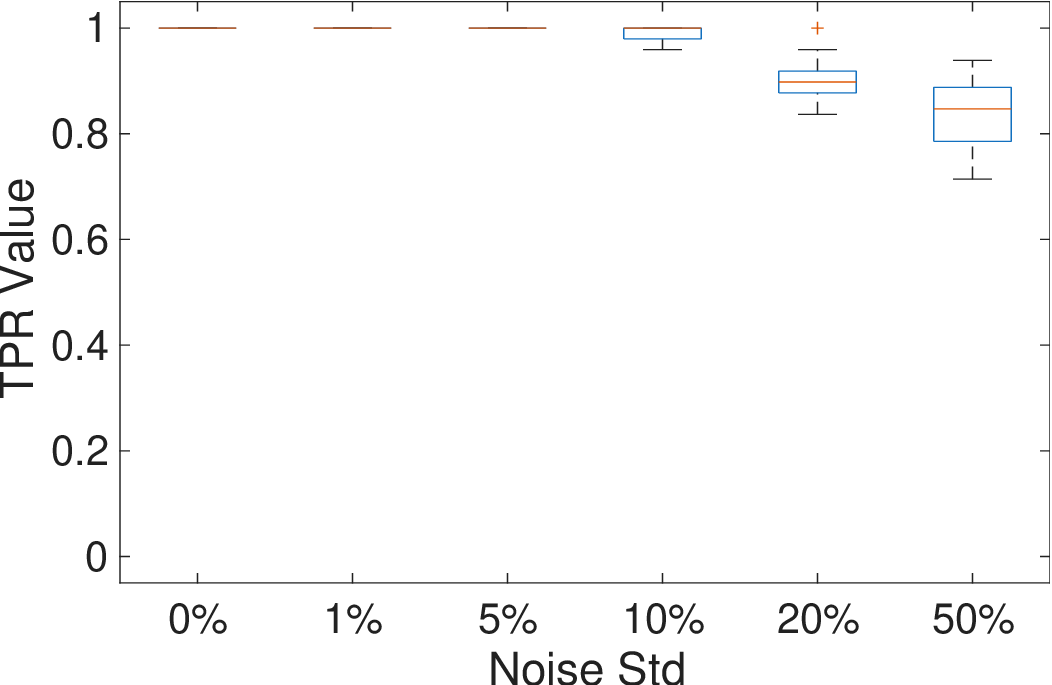}
  \caption{Prior-based (Weak form)}
  \label{subfig:no_prior_weak}
\end{subfigure}

\vspace{-0em}

\caption{
TPR results from twenty repeated identification trials for the three-body Hamiltonian system with different noise levels $\{0\%, 1\%, 5\%, 10\%, 20\%, 50\%\}$ under four configurations as mentioned above.
}

\label{fig:Hamiltoniansystem3body}
\end{figure}

\noindent\textbf{Identification results and performance.}  
We identify the dynamical system using both the baseline and prior-informed libraries across the four configurations. Under various noise levels, the TPR averaged over 20 trials is shown in Figure \ref{fig:Hamiltoniansystem3body}. With the baseline library, the 18 equations for position ($\dot{\mathbf{q}}_i$) and momentum ($\dot{\mathbf{p}}_i$) are identified independently from a comprehensive pool of 58 candidate features, necessitating a high-dimensional and unconstrained search. In contrast, the Hamiltonian prior couples the candidate terms into 12 structured groups via the energy basis $\mathcal{B}_H$, constraining the search to a single joint regression and drastically reducing the degrees of freedom. Across all noise levels, the proposed method exhibits stable performance. It correctly identifies the active support $\mathcal{S}^*_{\text{Ham}}$ corresponding to the 12 physical terms with a TPR of 1 up to 10\% noise, as shown in Figure \ref{fig:Hamiltoniansystem3body}(d). Conversely, results from the unconstrained baseline approach exhibit a noticeably reduced TPR across all noise levels.

\begin{figure}[ht!]
\centering

\begin{subfigure}{0.32\linewidth}
  \centering

  \includegraphics[trim={0.3cm 0.2cm 0.3cm 0.2cm}, clip, width=\linewidth]{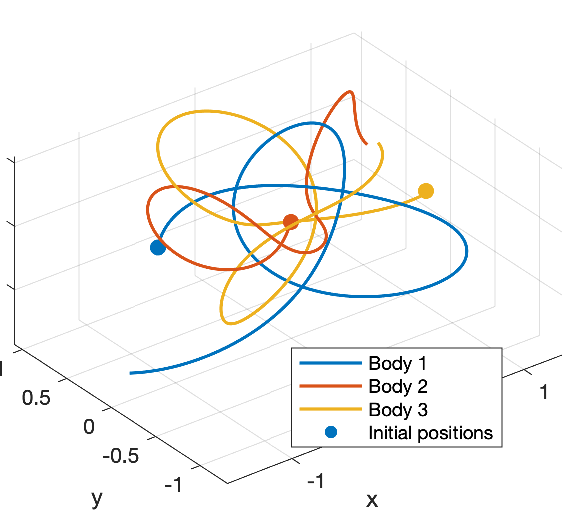}
  \caption{Given, 0\% noise}
  \label{subfig:true_00}
\end{subfigure}
\hfill
\begin{subfigure}{0.32\linewidth}
  \centering
  \includegraphics[trim={0.3cm 0.2cm 0.3cm 0.2cm}, clip, width=\linewidth]{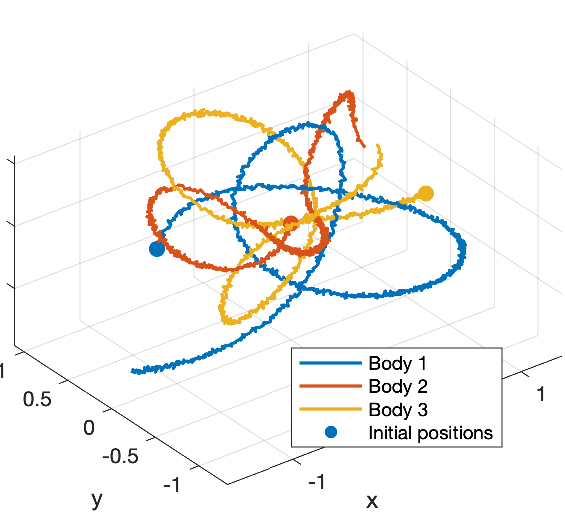}
  \caption{Given, 1\% noise}
  \label{subfig:true_10}
\end{subfigure}
\hfill
\begin{subfigure}{0.32\linewidth}
  \centering
  \includegraphics[trim={0.3cm 0.2cm 0.3cm 0.2cm}, clip, width=\linewidth]{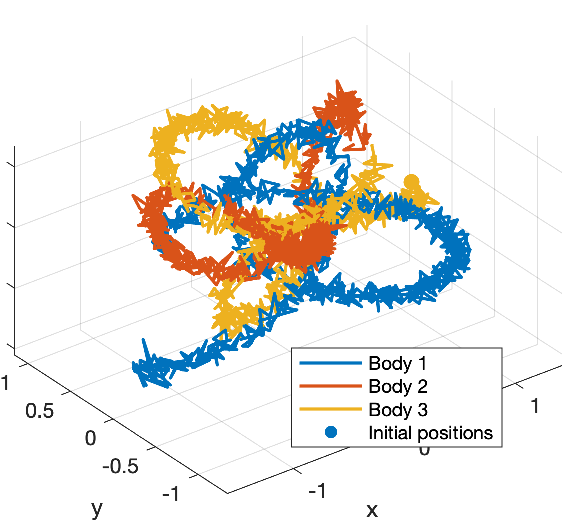}
  \caption{Given, 5\% noise}
  \label{subfig:true_50}
\end{subfigure}

\vspace{0.1em}

\begin{subfigure}{0.32\linewidth}
  \centering
  \includegraphics[trim={0.3cm 0.2cm 0.3cm 0.2cm}, clip, width=\linewidth]{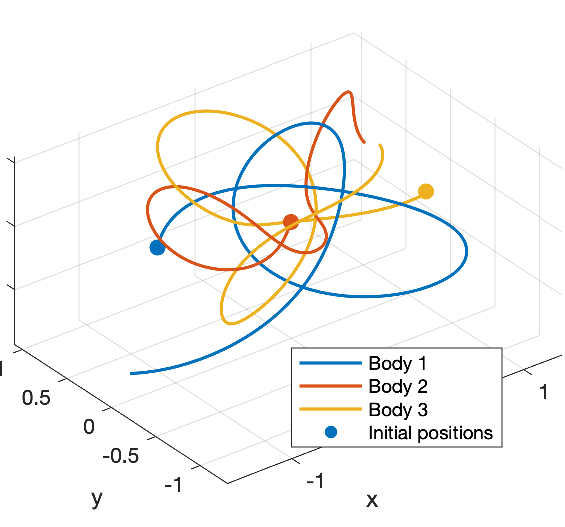}
  \caption{Identified, 0\% noise}
  \label{subfig:id_00}
\end{subfigure}
\hfill
\begin{subfigure}{0.32\linewidth}
  \centering
  \includegraphics[trim={0.3cm 0.2cm 0.3cm 0.2cm}, clip, width=\linewidth]{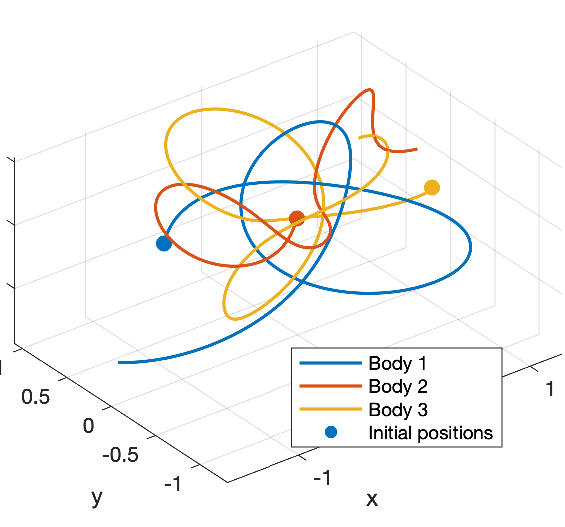}
  \caption{Identified, 1\% noise}
  \label{subfig:id_10}
\end{subfigure}
\hfill
\begin{subfigure}{0.32\linewidth}
  \centering
  \includegraphics[trim={0.3cm 0.2cm 0.3cm 0.2cm}, clip, width=\linewidth]{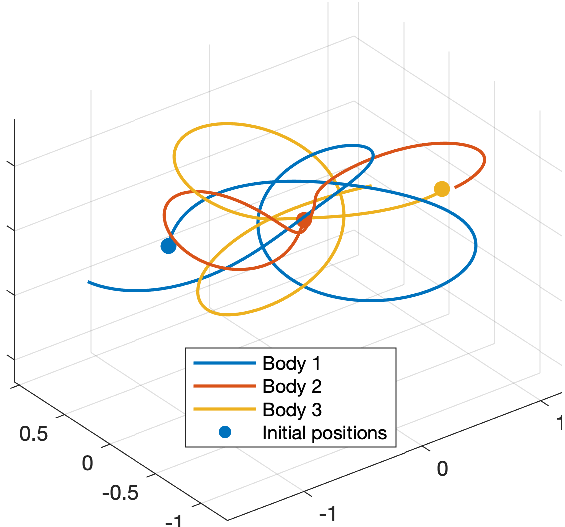}
  \caption{Identified, 5\% noise}
  \label{subfig:id_50}
\end{subfigure}

\caption{Comparison of trajectories of Three--Body Problem. Top row: True trajectories. Bottom row: Identified trajectories. Columns correspond to noise levels of 0\%, 1\%, and 5\% respectively.}
\label{fig:threebody_identified_trajectory_all_3D}
\end{figure}

Figure~\ref{fig:threebody_identified_trajectory_all_3D} compares the true and identified trajectories across noise levels (0\%, 1\%, and 5\%). While low-noise recoveries show excellent agreement, deviations emerge at 5\% noise. Notably, the framework successfully identifies the correct Hamiltonian structure (functional forms) even at this level; however, larger errors in coefficient estimation accumulate during integration, causing the observed spatial drift. Despite these quantitative deviations, the identified model maintains a qualitatively consistent orbital structure, demonstrating the robustness of the Hamiltonian prior.

\subsection{Conservation-Law prior}
\subsubsection{Inviscid Burgers’ equation} \label{numerical_exp:burgers}

We consider the inviscid Burgers' equation (see Section \ref{subsec:case_conservation}). Numerical data is generated by solving the conservative form \eqref{eq:Burgers_conservative} on a periodic domain $x \in [0, 1]$. The system is initialized with the mixed-wave condition $u_0(x) = 0.5(\sin(2\pi x) + \cos(2\pi x))$. Evolution is computed using a finite-volume scheme with the global Lax--Friedrichs flux to correctly capture shock formation. The spatial domain is discretized with $N_x = 500$ grid points, and trajectories are sampled over $N_t = 200$ time steps with $\Delta t = 0.001$.

\medskip
\noindent\textbf{Feature libraries.}  
To evaluate the efficacy of conservation-law priors, we construct two feature libraries and compare the baseline approach with the flux-constrained dictionary derived from the conservation form.

\noindent\textbf{Baseline (no-prior) dictionary.}  
The baseline dictionary $\mathcal{D}_{\text{base}}$ consists of polynomial derivative terms up to third order:
\begin{equation}\label{Burgers_dic_noprior}
\mathcal{D}_{\text{base}} = \big\{ u^k u_x,\, u^k u_{xx} \mid k=0,1,2,3 \big\} \cup \big\{ u^k \mid k=1,2,3 \big\}.
\end{equation}

This basis spans a wide space of nonlinear terms (e.g., $u\,u_x$, $u_{xx}$), offering sufficient expressiveness for general transport dynamics.  
However, lacking explicit conservation constraints, it may yield models with non-physical source terms.

\medskip
\noindent\textbf{Conservation-law prior dictionary.}  
Following the conservative structure identified in Section \ref{subsec:case_conservation}, we restrict admissible features to spatial derivatives of flux functions. Defining candidate fluxes $\mathcal{F} = \{ u, u^2, u^3 \}$, the flux-form dictionary is:
\begin{equation}\label{Burgers_flux_dic}
\mathcal{D}_{\text{flux}} = \big\{ (F)_x \mid F \in \mathcal{F} \big\} = \big\{ u_x, (u^2)_x, (u^3)_x \big\}.
\end{equation}
This construction guarantees that the recovered dynamics satisfy local flux balance by definition.

\begin{figure}[ht!]
\centering
\begin{subfigure}{0.45\linewidth}
  \includegraphics[width=\linewidth]{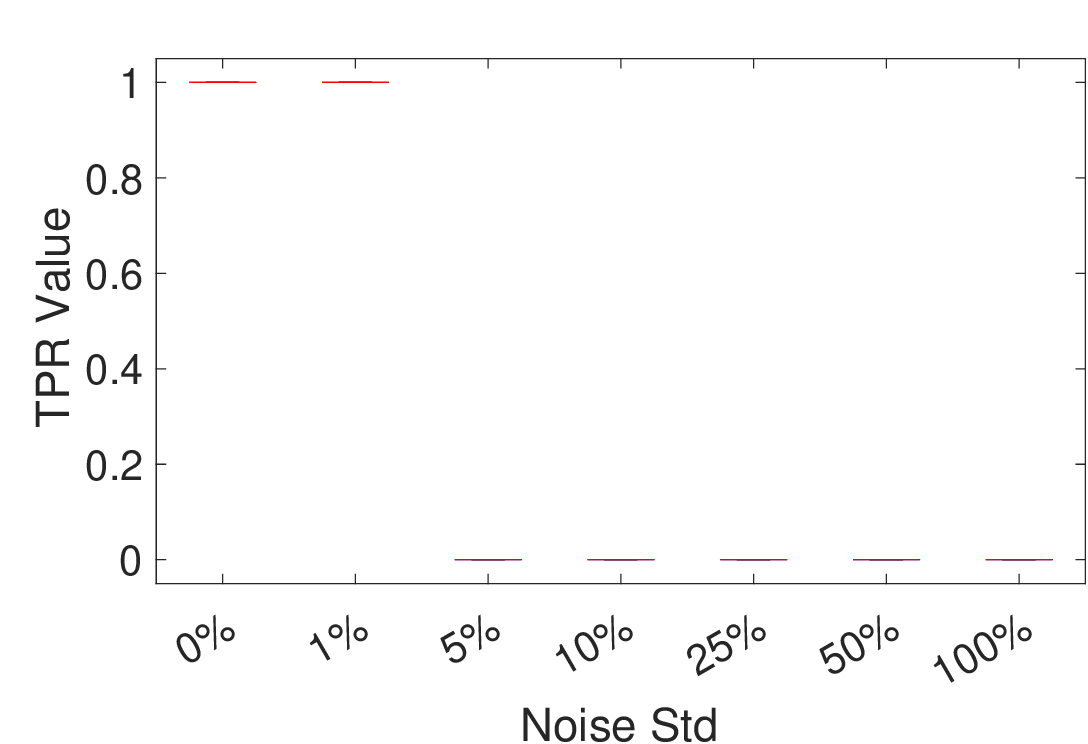}
  \caption{Separate (Strong form)}  
  \label{subfig:Burgers_SI_without_weak}
\end{subfigure}
\hfill  
\begin{subfigure}{0.45\linewidth}
  \includegraphics[width=\linewidth]{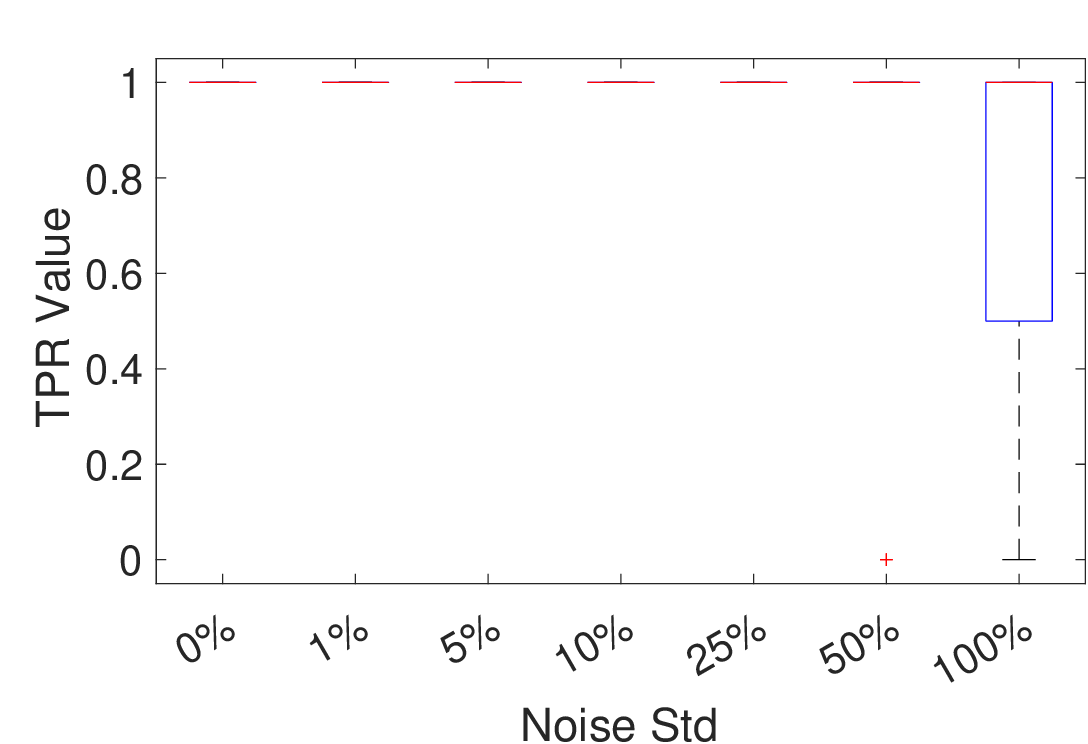}
  \caption{Separate (Weak form)}  
  \label{subfig:Burgers_SI_with_weak}
\end{subfigure}

\vspace{-0em}

\begin{subfigure}{0.45\linewidth}
  \includegraphics[width=\linewidth]{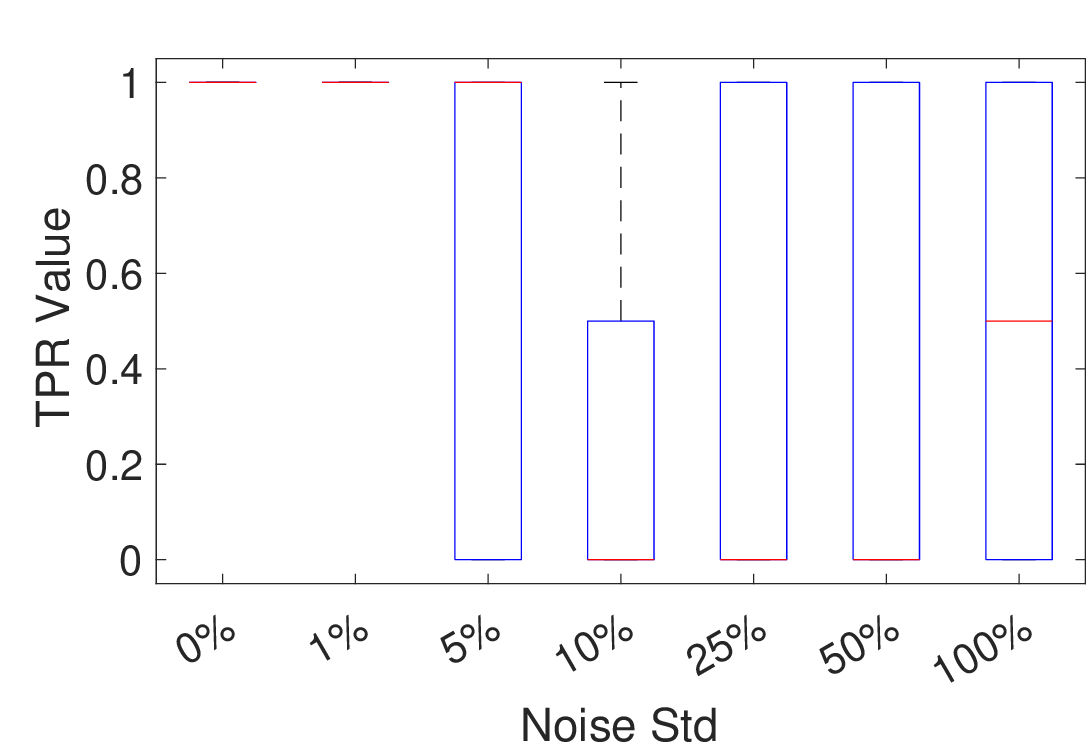}
  \caption{Prior-based (Strong form)}  
  \label{subfig:Burgers_PI_without_weak}
\end{subfigure}
\hfill
\begin{subfigure}{0.45\linewidth}
  \includegraphics[width=\linewidth]{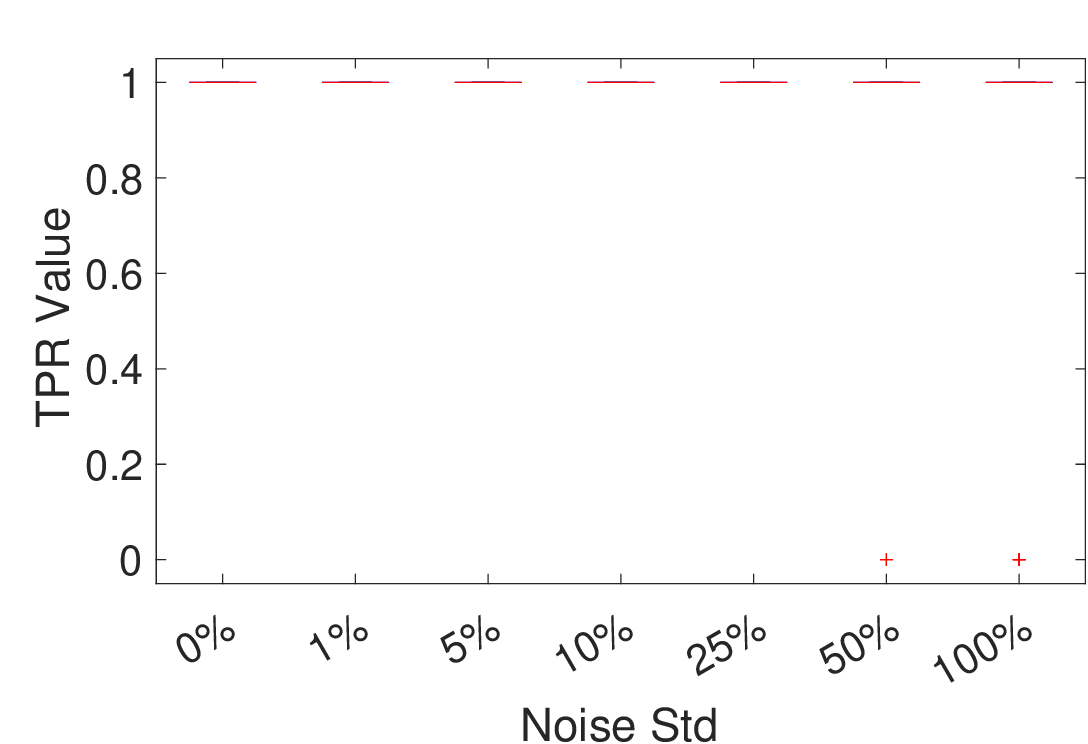}
  \caption{Prior-based (Weak form)}  
  \label{subfig:Burgers_PI_with_weak}
\end{subfigure}

\vspace{-0em}

\caption{TPR results from twenty repeated experiments for Burgers’ equation identification with different noise levels 
$\{0\%, 1\%, 5\%, 10\%, 25\%, 50\%, 100\%\}$ under four configurations as mentioned above.
}
\label{fig:burgers_identification}
\end{figure}

\noindent\textbf{Identification results and performance.}
We identify the governing PDE using both the baseline and flux-informed libraries across the four configurations. Under various noise levels, the TPR averaged over 20 trials is shown in Figure \ref{fig:burgers_identification}. With the baseline library, the regression must navigate a larger, unconstrained hypothesis space, increasing the risk of selecting non-conservative source terms or unphysical dissipation. In contrast, the conservation-law prior constrains the search strictly to flux gradients, significantly reducing the degrees of freedom. Across all noise levels, the proposed method exhibits stable performance. It correctly identifies the active support $\mathcal{S}^*_{\text{flux}}$ corresponding to the nonlinear flux component $(u^2)_x$ with a TPR of 1, as shown in Figure \ref{fig:burgers_identification}(d). Conversely, results with the baseline library exhibit a noticeably reduced TPR at high noise levels, often overfitting the noise with spurious terms. These results demonstrate the improved performance of our proposed method: the conservation-law prior ensures flux-balanced, physically consistent recovery of the nonlinear advection dynamics, while the weak-form integration substantially enhances numerical stability and robustness to noise.

\subsubsection{2D Shallow-Water equations} \label{numericalexp:2dSWE}

We simulate the conservative shallow-water equations defined in \eqref{eq:swe_flux}. The simulation is performed on a periodic domain $\Omega = [0, 1]^2$ with gravitational acceleration $g=9.81$. We initialize the system with a background water depth $H=1.5$. The initial condition is prescribed by a perturbed water height $h(x,y,0) = \max(H + \tilde{h}(x,y), 0.1H)$ and zero initial velocity $\mathbf{u} = \mathbf{0}$. The perturbation $\tilde{h}(x,y)$ consists of Gaussian anomalies and harmonic modes to ensure rich spectral content for identification:
\begin{equation}
\tilde{h}(x,y) = \sum_{k=1}^K A_k e^{-\alpha_k \|\mathbf{x}-\mathbf{\mu}_k\|^2} + 0.8\cos(2\pi x)\cos(2\pi y) + 0.3\cos(4\pi x)\cos(2\pi y),
\end{equation}
where $\mathbf{\mu}_k$, $A_k$, and $\alpha_k$ represent the centers, amplitudes, and concentration parameters of the Gaussian humps, respectively. In our primary numerical experiments, we employ $K=1$ with a centered Gaussian hump ($A_1=1.2, \alpha_1=800$) as shown in the reference setup.

This perturbation is centered and normalized to a standard deviation of $0.8H$. The velocity fields $(u_0, v_0)$ are initialized with sheared sinusoidal profiles balanced by pressure gradients:
\begin{equation}
\begin{aligned}
u_0 &= 0.8\sin(2\pi y) + 0.35\cos(2\pi x)\sin(2\pi y) + \beta\tfrac{g}{\omega}\partial_x h_0, \\
v_0 &= -0.6\sin(2\pi x) + 0.35\sin(2\pi x)\cos(2\pi y) + \beta\tfrac{g}{\omega}\partial_y h_0,
\end{aligned}
\end{equation}
where $\beta=0.40$ is the phase-alignment factor and $\omega = \sqrt{gH}|\mathbf{k}|$ represents the characteristic gravity-wave frequency. The system is evolved up to $T=0.3$ using a finite-volume scheme on a $100 \times 100$ grid with time step $\Delta t = 5 \times 10^{-4}$.

\medskip
\noindent\textbf{Feature libraries.}  To evaluate the efficacy of conservation-law priors, we construct two feature libraries and compare a standard polynomial-derivative library against a structure-preserving flux library.

\noindent\textbf{Baseline (no-prior) dictionary.}  
We employ a polynomial-derivative basis $\mathcal{D}_{\text{base}}$ constructed from products of state variables $\{1, h, u, v, u^2, v^2, uv, hu, hv, huv\}$ and their spatial derivatives $\{h_x, u_x, v_x, h_y, u_y, v_y\}$. This library captures convective and mixed nonlinearities (e.g., $u^2u_x$, $hv v_y$), representing general transport and coupling.  
However, this unconstrained basis offers no guarantees of local conservation and may admit non-physical interactions.

\medskip
\noindent\textbf{Conservation-law prior dictionary.}  
Based on the conservation form derived in Eq.~\eqref{eq:swe_flux}, we restrict the library to spatial divergences. Defining the candidate flux set $\mathcal{F}=\{h,hu,hv,u^2, \dots, uv^2\}$, the flux-form dictionary is:
\begin{equation}\label{SWE_dic_flux}
\mathcal{D}_{\text{flux}}
= \big\{ (F)_\alpha \mid F \in \mathcal{F},\ \alpha \in \{x,\,y\} \big\}.
\end{equation}
This prior enforces local conservation by construction while retaining expressiveness for advection and pressure coupling.

\begin{figure}[ht!]
\centering
\begin{subfigure}{0.45\linewidth}
  \includegraphics[width=\linewidth]{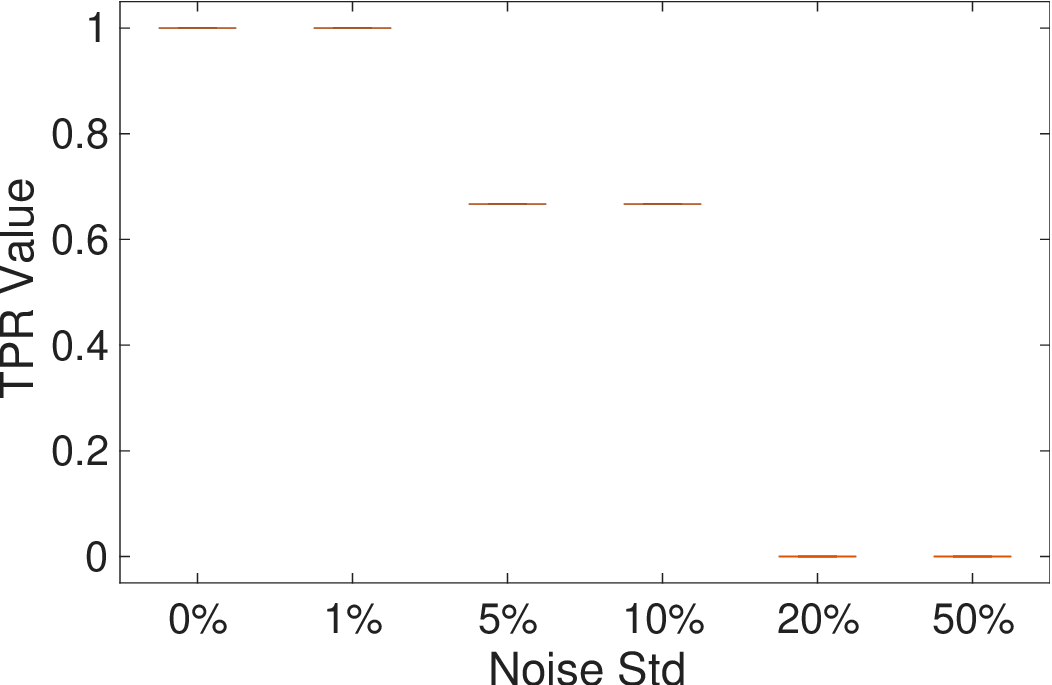}
  \caption{Separate (Strong form)}
  \label{subfig:SWE_SI_without_weak}
\end{subfigure}
\hfill
\begin{subfigure}{0.45\linewidth}
  \includegraphics[width=\linewidth]{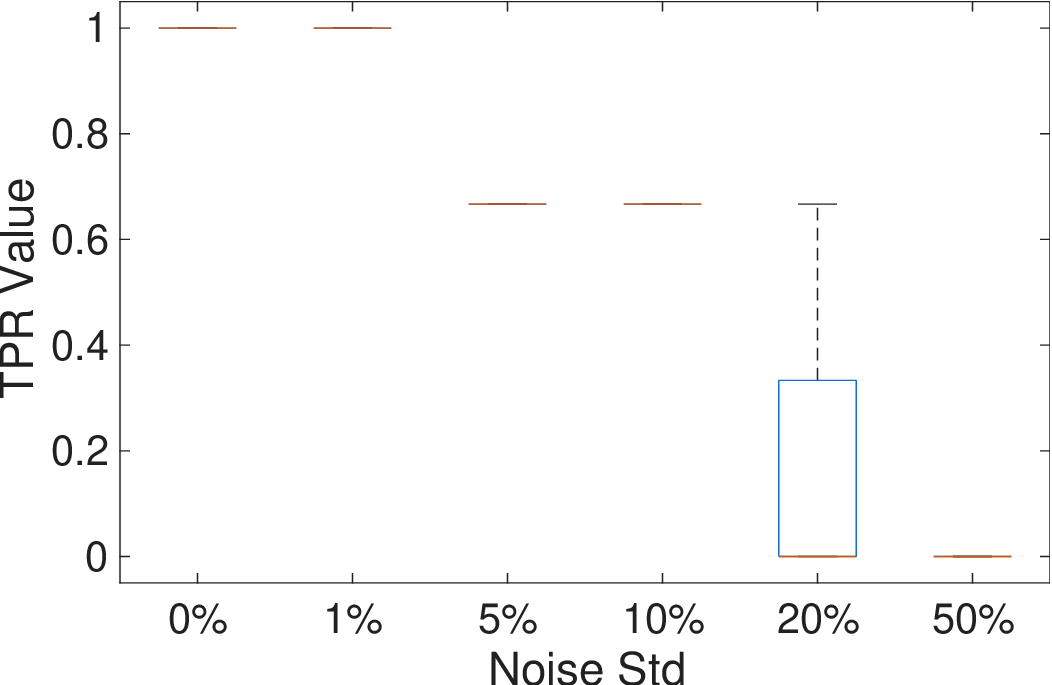}
  \caption{Separate (Weak form)}
  \label{subfig:SWE_SI_with_weak}
\end{subfigure}

\vspace{-0em}

\begin{subfigure}{0.45\linewidth}
  \includegraphics[width=\linewidth]{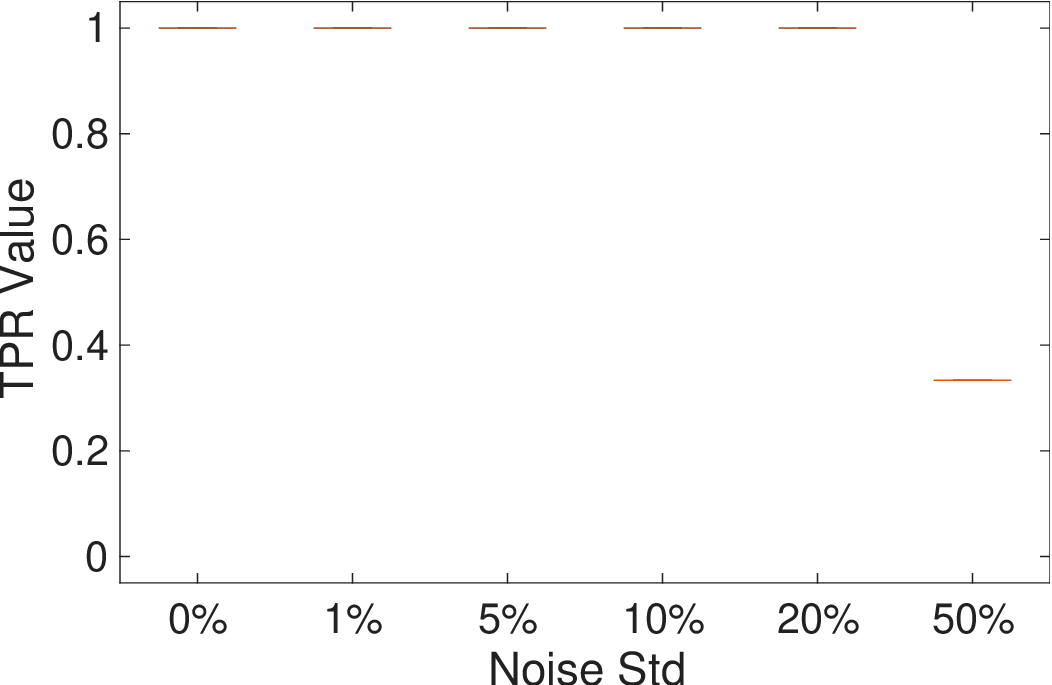}
  \caption{Prior-based (Strong form)}
  \label{subfig:SWE_PI_without_weak}
\end{subfigure}
\hfill
\begin{subfigure}{0.45\linewidth}
  \includegraphics[width=\linewidth]{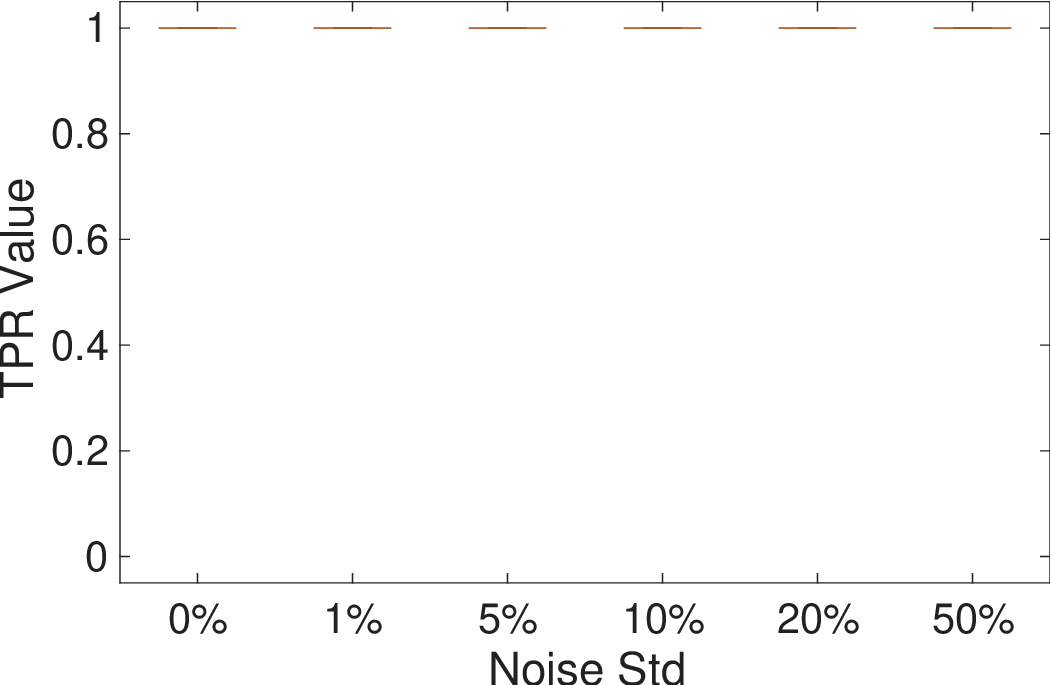}
  \caption{Prior-based (Weak form)}
  \label{subfig:SWE_PI_with_weak}
\end{subfigure}

\vspace{-0em}

\caption{TPR results from twenty repeated experiments for 2D SWE identification with different noise levels
$\{0\%, 1\%, 5\%, 10\%, 20\%, 50\%\}$ under four configurations as mentioned above. }
\label{fig:SWE_PI_with_weak_conservation}
\end{figure}

\begin{figure}[ht!]
\centering

\makebox[0.32\linewidth][c]{\small\textbf{$t=0$}}\hfill
\makebox[0.32\linewidth][c]{\small\textbf{$t=0.15$}}\hfill
\makebox[0.32\linewidth][c]{\small\textbf{$t=0.30$}}

\begin{subfigure}[t]{0.32\linewidth}
  \centering
  \includegraphics[width=\linewidth]{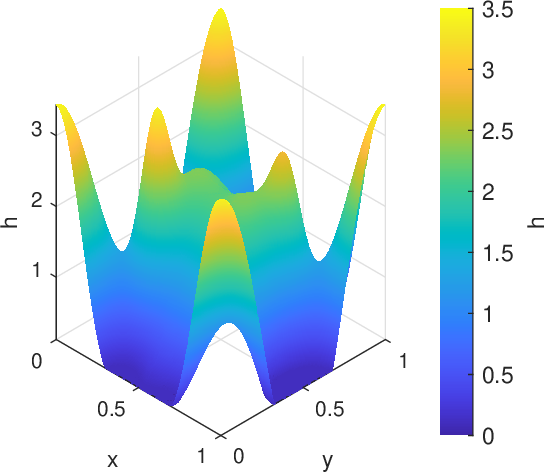}
  \label{subfig:swe_true_t0_00}
\end{subfigure}\hfill
\begin{subfigure}[t]{0.32\linewidth}
  \centering
  \includegraphics[width=\linewidth]{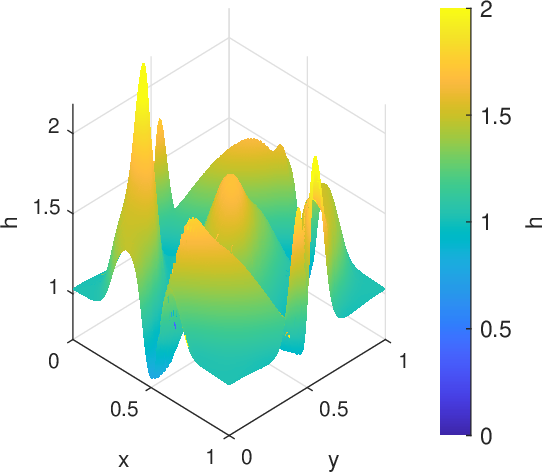}
  \label{subfig:swe_true_tmid_00}
\end{subfigure}\hfill
\begin{subfigure}[t]{0.32\linewidth}
  \centering
  \includegraphics[width=\linewidth]{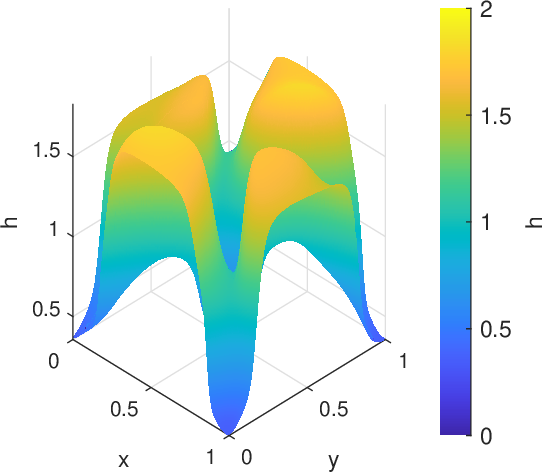}
  \label{subfig:swe_true_tend_00}
\end{subfigure}

\vspace{0.5em}

\begin{subfigure}[t]{0.32\linewidth}
  \centering
  \includegraphics[width=\linewidth]{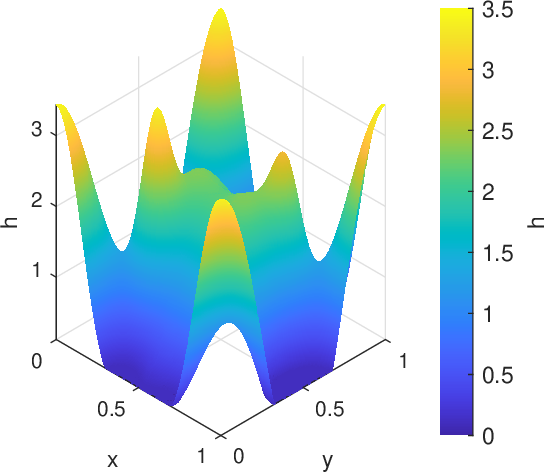}
  \label{subfig:swe_id_t0_00}
\end{subfigure}\hfill
\begin{subfigure}[t]{0.32\linewidth}
  \centering
  \includegraphics[width=\linewidth]{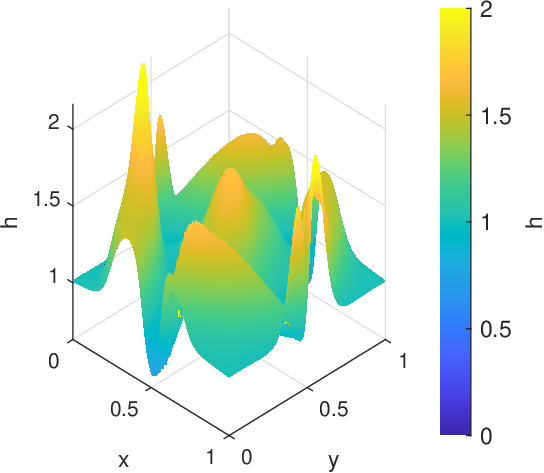}
  \label{subfig:swe_id_tmid_00}
\end{subfigure}\hfill
\begin{subfigure}[t]{0.32\linewidth}
  \centering
  \includegraphics[width=\linewidth]{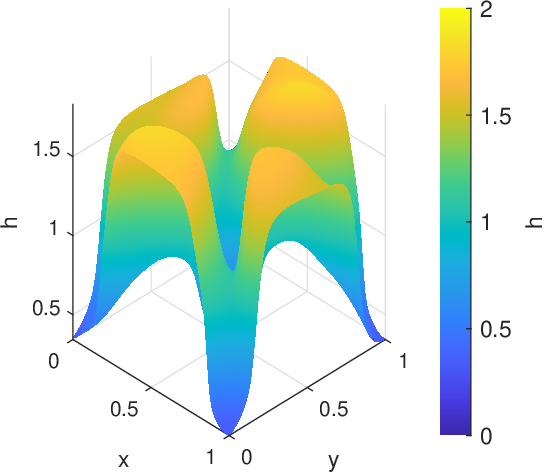}
  \label{subfig:swe_id_tend_00}
\end{subfigure}

\caption{Comparison of trajectories of 2D Shallow--Water equation under 0\% noise. Top row: Given trajectories. Bottom row: Identified trajectories. Columns correspond to different three time snapshots: $t=0$ (left column), mid (middle column), and end (right column).
}
\label{fig:SWE_h_trajectory_grid_00new}
\end{figure}

\begin{figure}[ht!]
\centering

\makebox[0.32\linewidth][c]{\small\textbf{$t=0$}}\hfill
\makebox[0.32\linewidth][c]{\small\textbf{$t=0.15$}}\hfill
\makebox[0.32\linewidth][c]{\small\textbf{$t=0.30$}}

\begin{subfigure}[t]{0.32\linewidth}
  \centering
  \includegraphics[width=\linewidth]{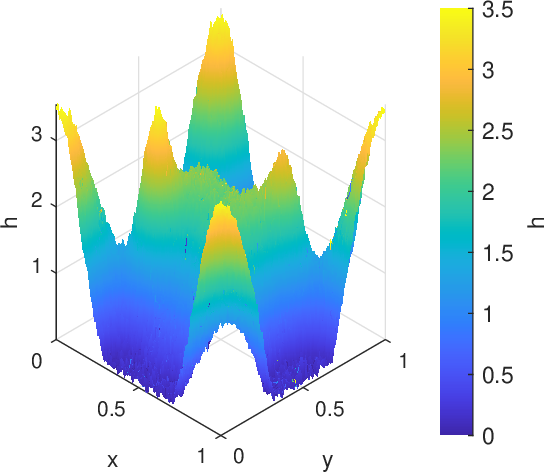}
  \label{subfig:swe_true_t0_10}
\end{subfigure}\hfill
\begin{subfigure}[t]{0.32\linewidth}
  \centering
  \includegraphics[width=\linewidth]{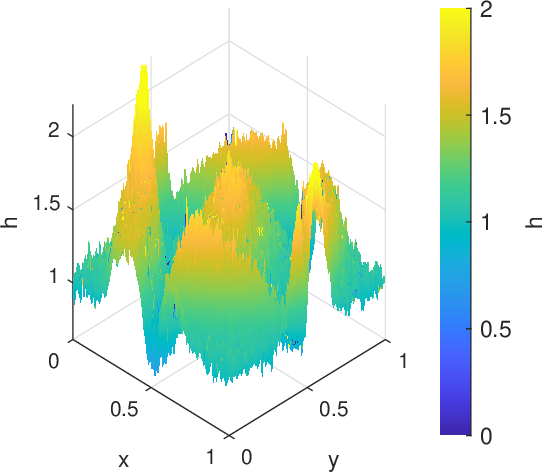}
  \label{subfig:swe_true_tmid_10}
\end{subfigure}\hfill
\begin{subfigure}[t]{0.32\linewidth}
  \centering
  \includegraphics[width=\linewidth]{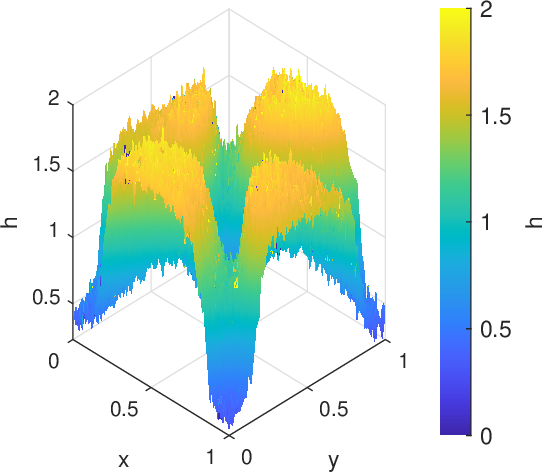}
  \label{subfig:swe_true_tend_10}
\end{subfigure}

\vspace{0.5em}

\begin{subfigure}[t]{0.32\linewidth}
  \centering
  \includegraphics[width=\linewidth]{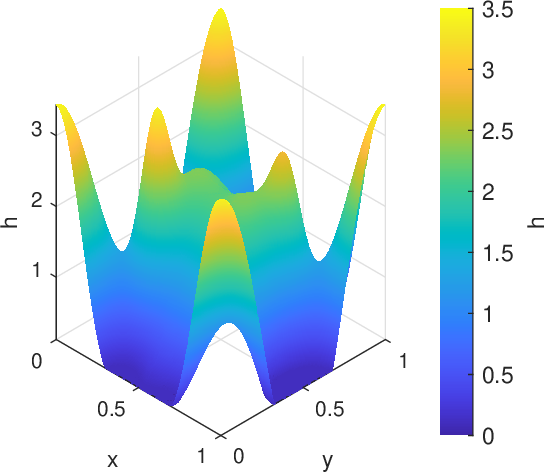}
  \label{subfig:swe_id_t0_10}
\end{subfigure}\hfill
\begin{subfigure}[t]{0.32\linewidth}
  \centering
  \includegraphics[width=\linewidth]{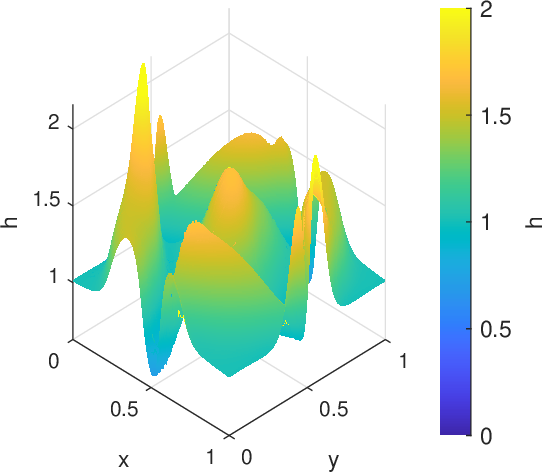}
  \label{subfig:swe_id_tmid_10}
\end{subfigure}\hfill
\begin{subfigure}[t]{0.32\linewidth}
  \centering
  \includegraphics[width=\linewidth]{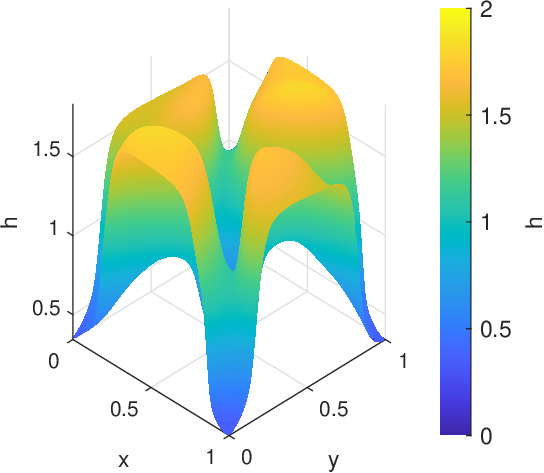}
  \label{subfig:swe_id_tend_10}
\end{subfigure}

\caption{Comparison of trajectories of 2D Shallow--Water equation under 10\% noise. Top row: Given trajectories. Bottom row: Identified trajectories. Columns correspond to different three time snapshots: $t=0$ (left column), $t=0.15$ (middle column), and $t=0.30$ (right column).
}
\label{fig:SWE_h_trajectory_grid_10new}
\end{figure}

\begin{figure}[ht!]
\centering

\makebox[0.32\linewidth][c]{\small\textbf{$t=0$}}\hfill
\makebox[0.32\linewidth][c]{\small\textbf{$t=0.15$}}\hfill
\makebox[0.32\linewidth][c]{\small\textbf{$t=0.30$}}

\begin{subfigure}[t]{0.32\linewidth}
  \centering
  \includegraphics[width=\linewidth]{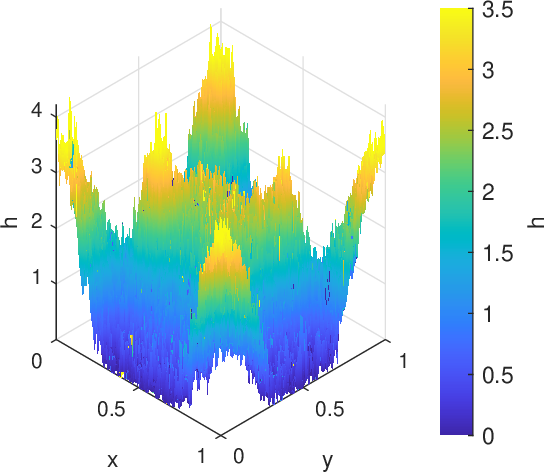}
  \label{subfig:swe_true_t0_50}
\end{subfigure}\hfill
\begin{subfigure}[t]{0.32\linewidth}
  \centering
  \includegraphics[width=\linewidth]{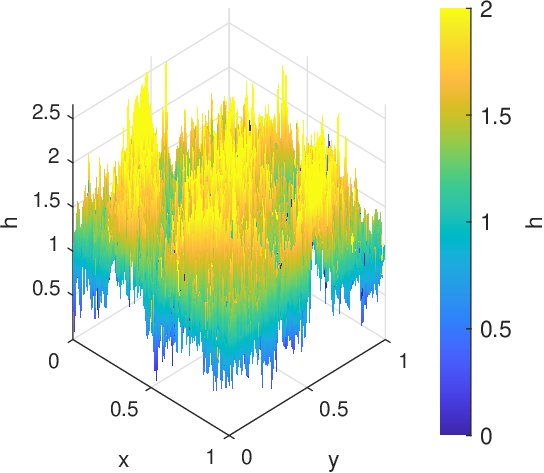}
  \label{subfig:swe_true_tmid_50}
\end{subfigure}\hfill
\begin{subfigure}[t]{0.32\linewidth}
  \centering
  \includegraphics[width=\linewidth]{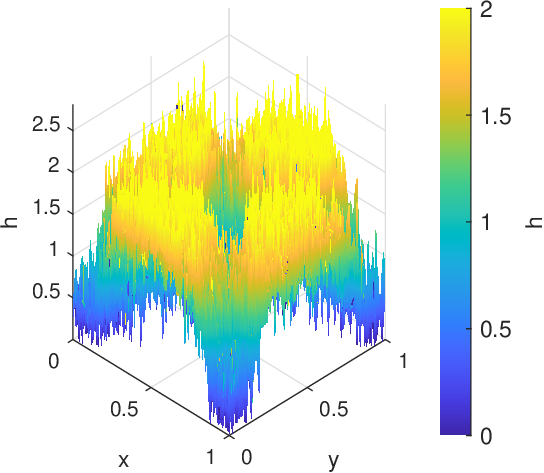}
  \label{subfig:swe_true_tend_50}
\end{subfigure}

\vspace{0.5em}

\begin{subfigure}[t]{0.32\linewidth}
  \centering
  \includegraphics[width=\linewidth]{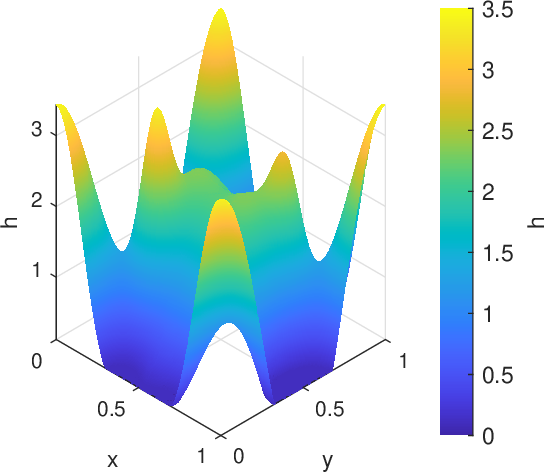}
  \label{subfig:swe_id_t0_50}
\end{subfigure}\hfill
\begin{subfigure}[t]{0.32\linewidth}
  \centering
  \includegraphics[width=\linewidth]{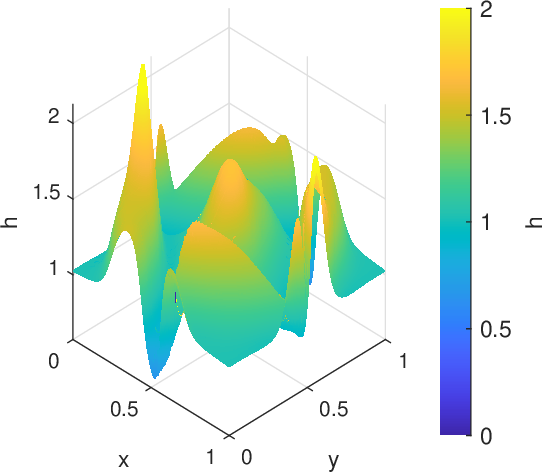}
  \label{subfig:swe_id_tmid_50}
\end{subfigure}\hfill
\begin{subfigure}[t]{0.32\linewidth}
  \centering
  \includegraphics[width=\linewidth]{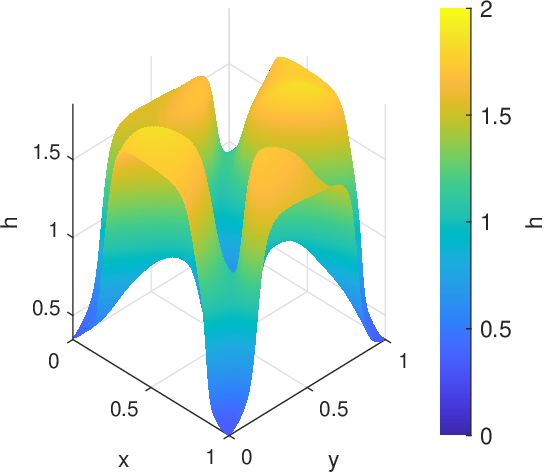}
  \label{subfig:swe_id_tend_50}
\end{subfigure}

\caption{Comparison of trajectories of 2D Shallow--Water equation under 50\% noise. Top row: Given trajectories. Bottom row: Identified trajectories. Columns correspond to different three time snapshots: $t=0$ (left column), $t=0.15$ (middle column), and $t=0.30$ (right column).}
\label{fig:SWE_h_trajectory_grid_50new}
\end{figure}

\noindent\textbf{Identification results and performance.}
We identify the governing PDE system using both the baseline and flux-informed libraries across the four configurations. Under various noise levels, the TPR averaged over 20 trials is shown in Figure \ref{fig:SWE_PI_with_weak_conservation}. With the baseline library, the equations for $h$, $hu$, and $hv$ are identified independently from a large pool of 60 candidate features. This unconstrained approach requires selecting component-wise sparsities ($\theta_h^*=2$, $\theta_{hu}^*=3$, and $\theta_{hv}^*=3$) by using RR and significantly increases the risk of selecting spurious, non-physical terms. In contrast, the conservation-law prior restricts the search strictly to spatial divergences via $\mathcal{D}_{\text{flux}}$, reducing the degrees of freedom through a single joint regression. Across all noise levels, the proposed method exhibits stable performance. It correctly identifies the active support $\mathcal{S}^*_{\text{flux}}$ corresponding to the 8 conservative physical terms with a TPR of 1 across all noise levels. Notably, it not only recovers advective fluxes with coefficients closely approximating unity, but also accurately identifies the hydrostatic pressure terms $(h^2)_x$ and $(h^2)_y$ with coefficients matching the physical constant $\frac{1}{2}g$. Conversely, the unconstrained baseline approach exhibits a noticeably reduced TPR at high noise levels. These results demonstrate the improved performance of our proposed method: the conservation-law prior ensures structure-preserving, physically consistent recovery, while the weak-form integration substantially enhances numerical stability and robustness to noise.  Furthermore, as illustrated in Figures \ref{fig:SWE_h_trajectory_grid_00new}--\ref{fig:SWE_h_trajectory_grid_50new}, the identified models closely reproduce the spatiotemporal evolution of the water height field $h(x,y,t)$ at representative snapshots ($t \in \{0, 0.15, 0.3\}$), preserving conservative dynamics even under severe noise up to 50\%.

We quantify the reconstruction accuracy using two metrics: the relative coefficient error $E_{\mathbf{c}} = \|\mathbf{c}^* - \mathbf{c}_{\text{true}}\|_2 / \|\mathbf{c}_{\text{true}}\|_2$, where $\mathbf{c}^*$ and $\mathbf{c}_{\text{true}}$ denote the identified and true coefficient vectors, and the total relative $\ell_2$ error for the state variables $
    E_{\text{total}} = \|\mathbf{q}_{\text{id}} - \mathbf{q}_{\text{true}}\|_2/\|\mathbf{q}_{\text{true}}\|_2$,
where $\mathbf{q} = [h, u, v]^\top$ represents the combined state vector. The resulting errors are summarized in Table~\ref{tab:SWE_E2_noise}. Notably, our method achieves consistently low errors even under substantial measurement noise; as shown in the table, the total relative error remains below $3\%$ (specifically $2.995 \times 10^{-2}$) even at a $50\%$ noise level, demonstrating the exceptional robustness of the proposed framework.

\begin{table}[ht!]
\centering
\caption{Relative $\ell^2$ errors of the identified 2D shallow-water equations (SWE) solutions under different noise levels.}
\label{tab:SWE_E2_noise}
\begin{tabular}{ccccc}
\toprule
\textbf{Noise Level (\%)} & \textbf{h} & \textbf{u} & \textbf{v} & \textbf{Total} \\
\midrule
0  & $6.903\times10^{-3}$ & $2.017\times10^{-2}$ & $2.477\times10^{-2}$ & $1.748\times10^{-2}$ \\
10 & $6.933\times10^{-3}$ & $2.039\times10^{-2}$ & $2.469\times10^{-2}$ & $1.753\times10^{-2}$ \\
50 & $1.338\times10^{-2}$ & $3.788\times10^{-2}$ & $3.801\times10^{-2}$ & $2.995\times10^{-2}$ \\
\bottomrule
\end{tabular}
\end{table}

\subsection{Energy-Dissipation Prior}

\subsubsection{Diffusion Equation}\label{numericalexp:diffu}

We consider the one-dimensional diffusion equation as defined in Section~\ref{subsec:case_energy_diffusioneq}(a) 1D Diffusion Equation. We generate numerical data by solving ~\eqref{eq:diffusion_strong} on the domain $[0, 1]$ with periodic boundary conditions. The diffusivity is set to $\nu=0.02$. The initial condition is a composite profile $u(x,0) = e^{-600(x-0.5)^2} + 0.2\sin(4\pi x)$, combining a localized Gaussian with a harmonic perturbation. Time integration is performed over the interval $t \in [0, 0.2]$ using the Forward-Time Central-Space (FTCS) explicit scheme. The spatial domain is discretized with $N_x = 500$ grid points ($\Delta x = 0.002$). To ensure numerical stability, the time step is strictly set to $\Delta t = 2.5 \times 10^{-5}$. 

\medskip

\noindent\textbf{Feature Libraries.}  
We construct two feature libraries to compare the baseline approach with the energy-informed prior defined in Section \ref{subsec:case_energy_diffusioneq}.

\noindent\textbf{Baseline (no-prior) dictionary.}  
The unconstrained library $\mathcal{D}_{\text{base}}$ contains state and derivative terms up to the second order:
$$\mathcal{D}_{\text{base}}
= \{\, u,\, u^2,\, u_x,\, u_x^2,\, u_{xx},\, (u_{xx})^2 \,\}.$$
While this library accommodates linear diffusion $(u_{xx})$, it also includes nonlinear terms (e.g., $u_x^2, (u_{xx})^2$) that may not correspond to valid gradient flow dynamics derived from a well-posed energy functional.

\noindent\textbf{Energy–Dissipation prior dictionary.}  
Following the gradient flow structure derived in Section \ref{subsec:case_energy_diffusioneq}, we construct the valid library by identifying variational derivatives of admissible energy densities. We define the candidate energy library $\mathcal{D}_{\text{energy}}^{\text{valid}}$ using convex terms consistent with the Dirichlet energy:
\begin{equation}\label{diff_dic_valid}
\mathcal{D}_{\text{energy}}^{\text{valid}}
= \{\, u^2,\, u_x^2,\, (u_{xx})^2 \,\}\,
\end{equation}
The corresponding PDE-level library $\mathcal{D}_{\text{GF}}^{\text{valid}}$ is obtained via the variational mapping: 

\begin{equation}\label{diff_gf_valid}
\mathcal{D}_{\text{GF}}^{\text{valid}}
= \{\, -u,\, +u_{xx},\, -u_{xxxx} \,\}.
\end{equation}

This framework precludes non-variational and indefinite terms, restricting regression to a subset of thermodynamically consistent operators. This refinement ensures physical interpretability and numerical stability by construction.

\begin{figure}[t!]
\centering
\begin{subfigure}{0.45\linewidth}
  \includegraphics[width=\linewidth]{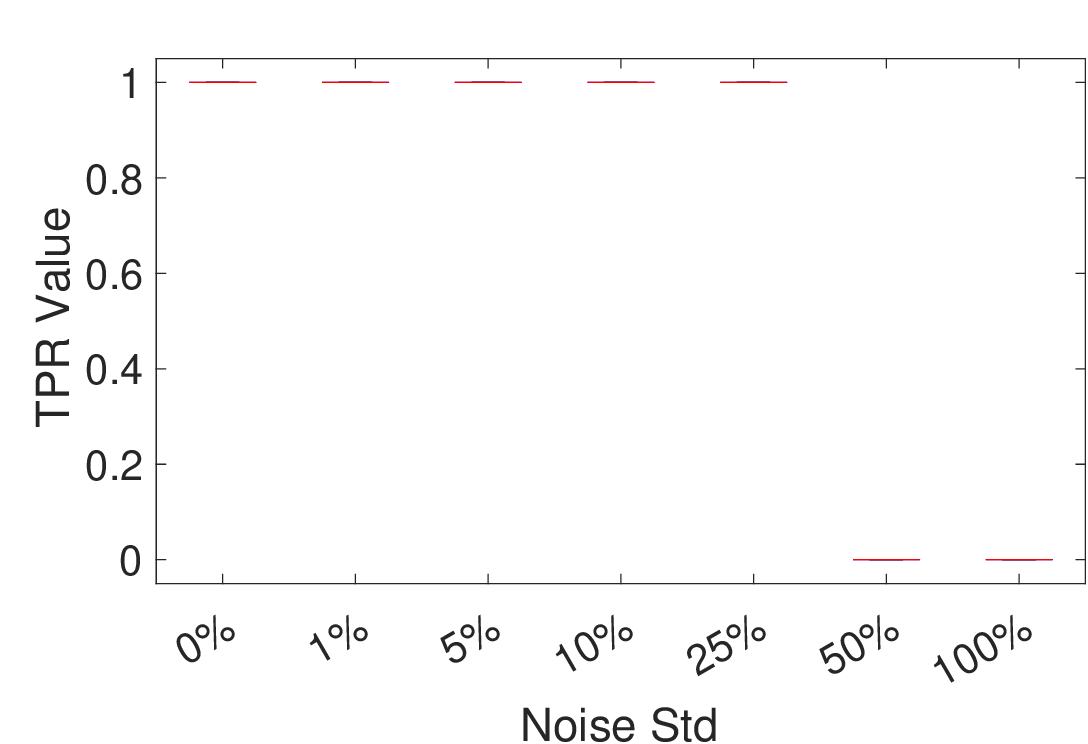}
  \caption{Separate (Strong form)}  
  \label{subfig:diff_SI_without_weak}
\end{subfigure}
\hfill  
\begin{subfigure}{0.45\linewidth}
  \includegraphics[width=\linewidth]{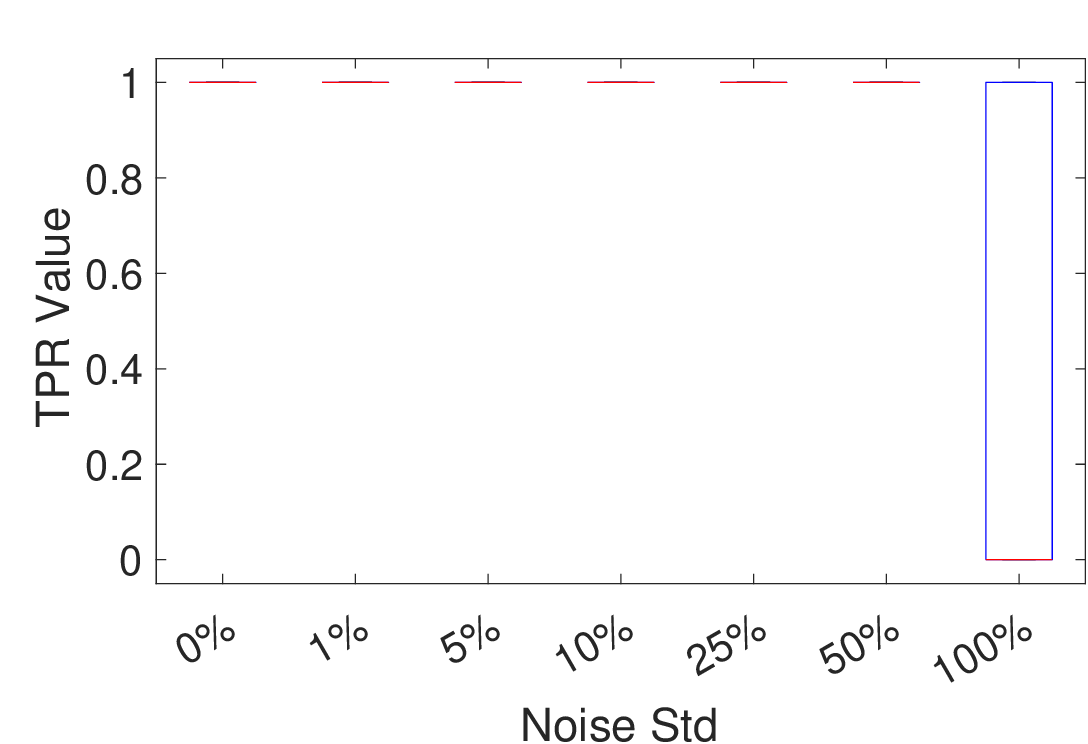}
  \caption{Separate (Weak form)}  
  \label{subfig:diff_SI_with_weak}
\end{subfigure}

\vspace{-0em}

\begin{subfigure}{0.45\linewidth}
  \includegraphics[width=\linewidth]{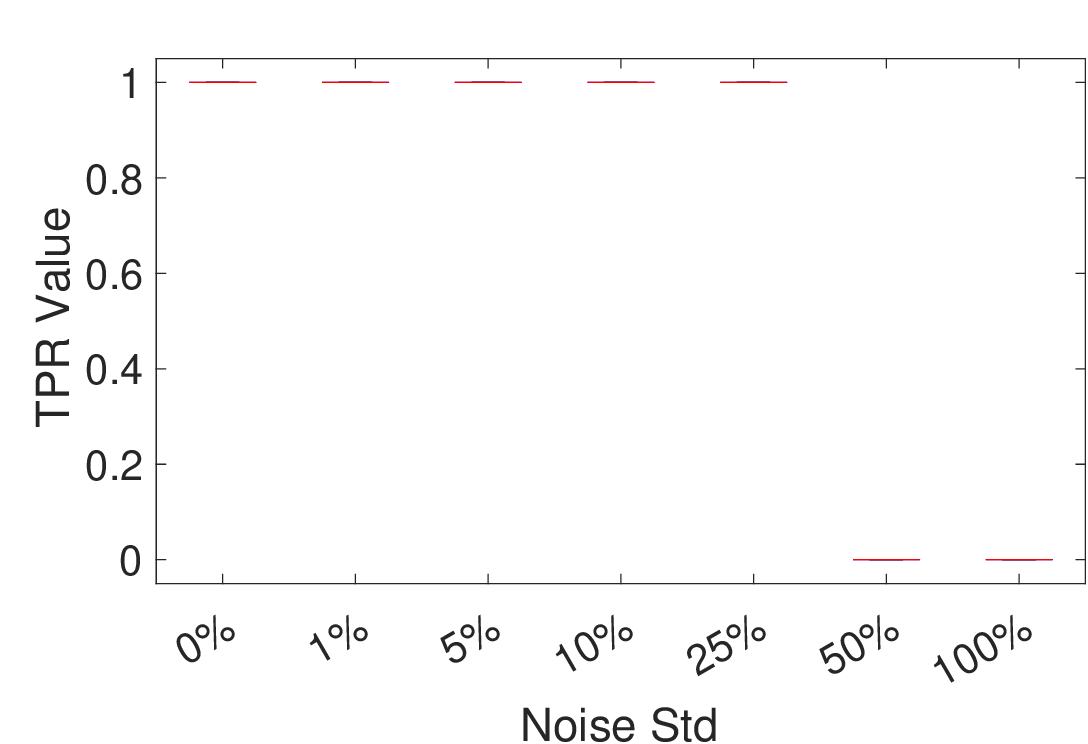}
  \caption{Prior-based (Strong form)}  
  \label{subfig:diff_PI_without_weak}
\end{subfigure}
\hfill
\begin{subfigure}{0.45\linewidth}
  \includegraphics[width=\linewidth]{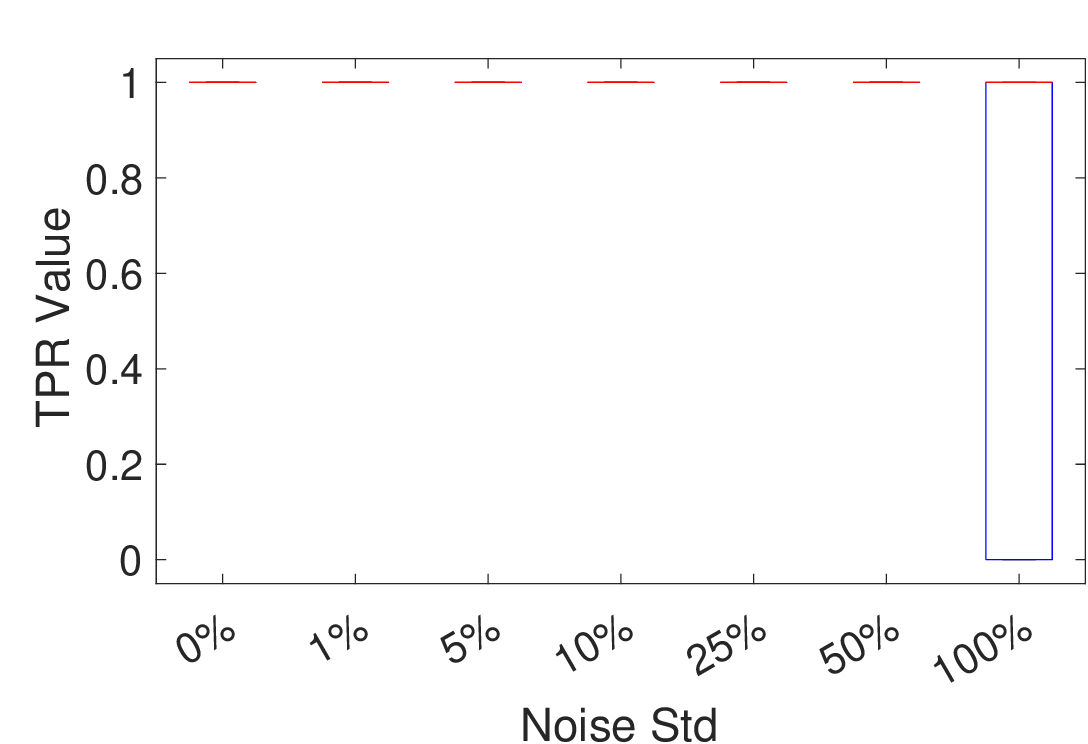}
  \caption{Prior-based (Weak form)}  
  \label{subfig:diff_PI_with_weak}
\end{subfigure}

\vspace{-0em}

\caption{TPR results from twenty repeated experiments for diffusion equation identification with varying noise levels 
$\{0\%, 1\%, 5\%, 10\%, 25\%, 50\%, 100\%\}$ under four configurations as mentioned above. }
\label{fig:diffusion_identification}
\end{figure}

\noindent\textbf{Identification results and performance.}
We identify the governing PDE using both the baseline and energy-informed libraries across the four configurations. Under various noise levels, the TPR averaged over 20 trials is shown in Figure \ref{fig:diffusion_identification}. With the baseline library, the regression must navigate an unconstrained hypothesis space, increasing the risk of selecting non-variational or indefinite terms that violate fundamental thermodynamic principles. In contrast, the Energy--Dissipation prior restricts the search strictly to valid gradient flow dynamics via $\mathcal{D}_{\text{GF}}^{\text{valid}}$, reducing the degrees of freedom. It correctly identifies the active support corresponding to the Laplacian operator $u_{xx}$ and recovers the model $u_t = \nu^\star u_{xx}$ with the identified coefficient $\nu^\star$ closely approximating the ground truth $\nu=0.02$, which inherently guarantees monotonic energy decay ($\frac{dE}{dt} \le 0$). When comparing configuration (d) to the baseline in (b), both weak-form approaches demonstrate strong robustness, maintaining perfect identification up to 50\% noise. However, at the extreme 100\% noise level, the advantage of the proposed method becomes statistically apparent: while the unconstrained baseline in (b) exhibits high variance and a dropped median (the red segment) TPR, the prior-informed method in (d) manages to maintain a higher median recovery rate (median TPR=1). This subtle but important difference indicates that under severe noise conditions, the thermodynamic constraints provide a crucial stabilizing effect, preventing the selection of physically invalid terms that unconstrained methods might otherwise fit to noise.

\subsubsection{Allen--Cahn Equation}\label{numericalexp:allen_cahn}

For the Allen-Cahn system described in Section \ref{subsec:case_energy_allencahn}(b), we simulate the dynamics on a periodic domain $[0, 2\pi)$ using a spectral method. The system is initialized with a multi-mode perturbation $u_0(x) = \cos x + \cos(2x) + 0.5\cos(3x)$ and evolved up to $T=2.0$ with a time step of $10^{-3}$.

\noindent\textbf{Feature Libraries.}
\label{sec:allen_cahn_feature_library}
We compare an unconstrained library against the energy-based construction detailed in Section \ref{subsec:case_energy} (b) 1D Allen–Cahn Equation.

\noindent\textbf{Baseline (no-prior) dictionary.}
We define an unconstrained library $\mathcal{D}_{\text{base}}$ containing $M=18$ candidate terms. The library is constructed from polynomial powers of $u$ and its spatial derivatives:
\begin{equation}
    \mathcal{D}_{\text{base}} = \big\{ u^p, (u^p)_x, (u^p)_{xx} \big\}_{p=1}^4 \cup \big\{ (u_x)^p, (u_{xx})^p \big\}_{p=2}^4,
\end{equation}
While this overcomplete library captures both reaction and diffusion dynamics, it may admit non-physical or anti-dissipative terms, leading to ill-posed models.

\noindent\textbf{Energy–Dissipation prior dictionary.}
Consistent with the gradient flow structure discussed in Section \ref{subsec:case_energy_allencahn}, we construct the candidate library by identifying variational derivatives of admissible energy densities. We define the valid energy library $\mathcal{D}_{\text{energy}}^{\text{valid}}$ using even-powered polynomials and squared derivatives:
\begin{equation}\label{ac_dic_valid}
\mathcal{D}_{\text{energy}}^{\text{valid}} = \big\{ u^2,\, u^4,\, (u_x)^2,\, (u_x)^4,\, (u_{xx})^2,\, (u_{xx})^4 \big\}.
\end{equation}

The corresponding gradient flow library $\mathcal{D}_{\text{GF}}^{\text{valid}}$ is obtained via the variational derivative as shown following: 

\begin{equation}\label{ac_gf_valid}
\begin{split}
\mathcal{D}_{\text{GF}}^{\text{valid}} = \big\{ & -2u,\, -4u^3,\, +2u_{xx},\, +12u_x^2u_{xx}, \\
& -2u_{xxxx},\, -(24u_{xx}u_{xxx}^2 + 12u_{xx}^2u_{xxxx}) \big\}.
\end{split}
\end{equation}

This construction restricts the search space to a thermodynamically consistent and physically interpretable model class.

\begin{figure}[ht!]
\centering
\begin{subfigure}{0.45\linewidth}
  \includegraphics[width=\linewidth]{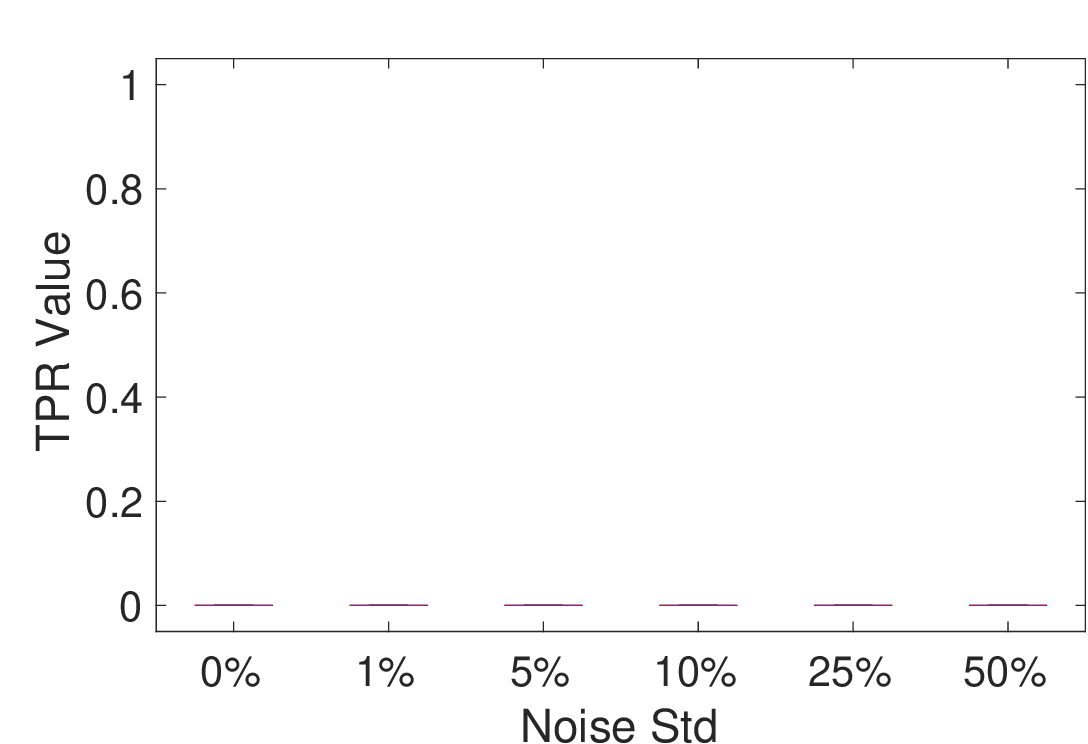}
  \caption{Separate (Strong form)}  
  \label{subfig:AC_SI_Strong}
\end{subfigure}
\hfill  
\begin{subfigure}{0.45\linewidth}
  \includegraphics[width=\linewidth]{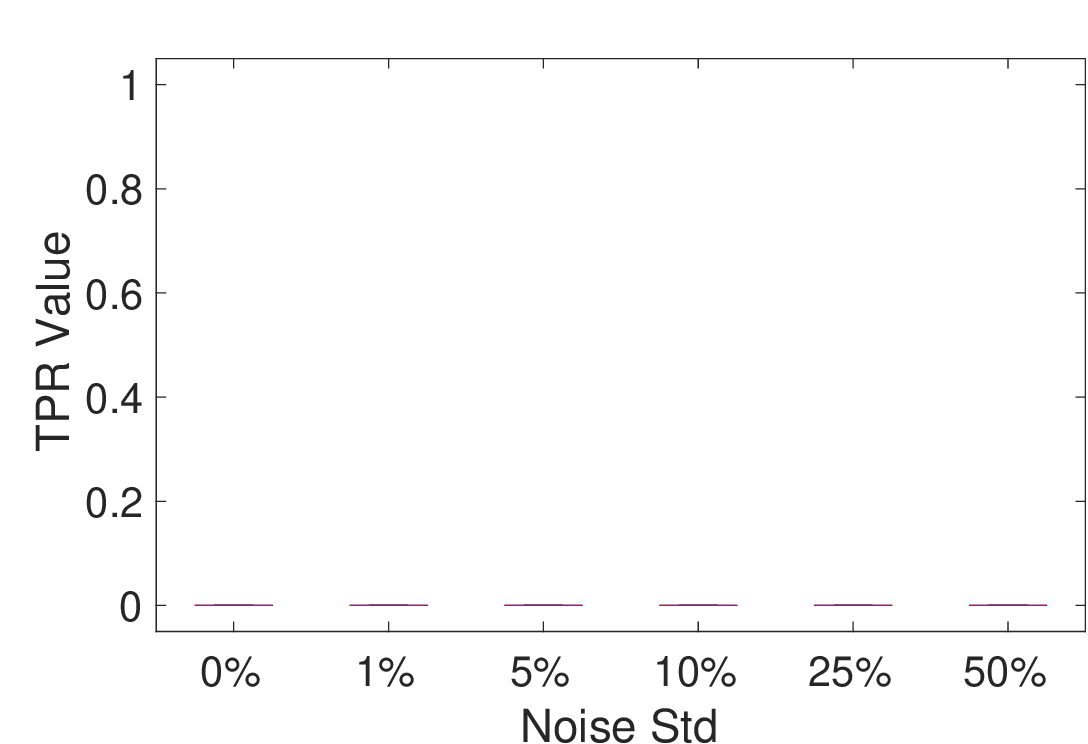}
  \caption{Separate (Weak form)}  
  \label{subfig:AC_SI_Weak}
\end{subfigure}

\vspace{-0em}

\begin{subfigure}{0.45\linewidth}
  \includegraphics[width=\linewidth]{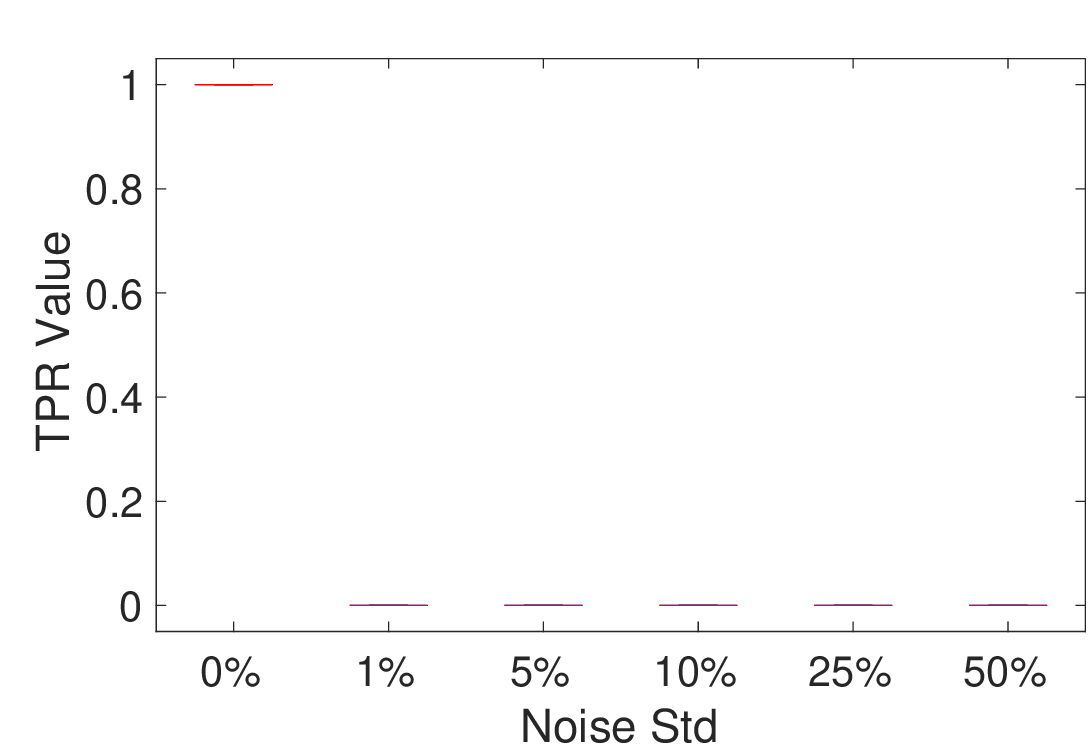}
  \caption{Prior-based (Strong form)}  
  \label{subfig:AC_PI_Strong}
\end{subfigure}
\hfill
\begin{subfigure}{0.45\linewidth}
  \includegraphics[width=\linewidth]{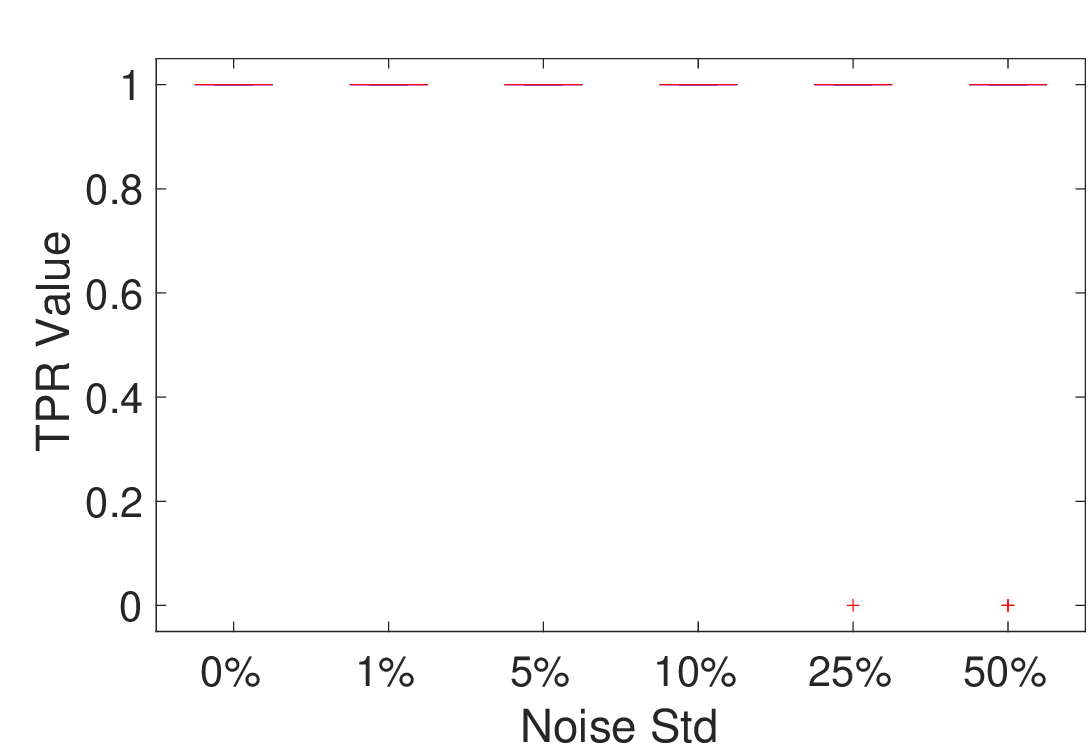}
  \caption{Prior-based (Weak form)}  
  \label{subfig:AC_PI_Weak}
\end{subfigure}

\vspace{-0em}

\caption{TPR results from twenty repeated experiments for Allen--Cahn equation identification with different noise levels $\{0\%, 1\%, 5\%, 10\%, 20\%, 50\% \}$ under four configurations as mentioned above. }
\label{fig:AC_identification}
\end{figure}

\begin{figure}[ht!]
\centering

\makebox[0.32\linewidth][c]{\small 0\% Noise}\hfill
\makebox[0.32\linewidth][c]{\small 5\% Noise}\hfill
\makebox[0.32\linewidth][c]{\small 10\% Noise}

\vspace{0.2em}

\begin{subfigure}[t]{0.32\linewidth}
  \centering

  \includegraphics[width=\linewidth]{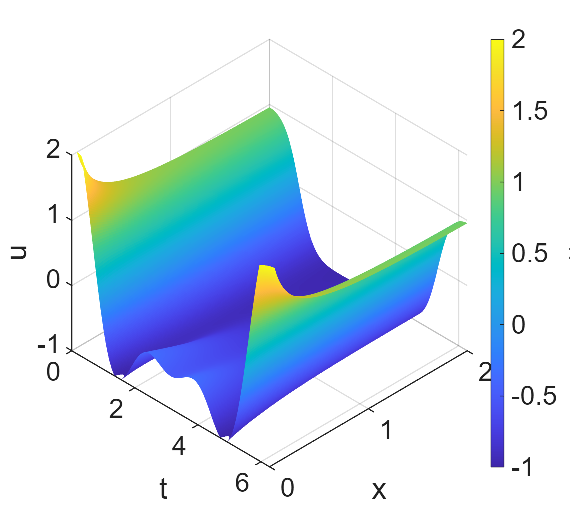}
  \label{subfig:AC_noise00_true}
\end{subfigure}\hfill
\begin{subfigure}[t]{0.32\linewidth}
  \centering

  \includegraphics[width=\linewidth]{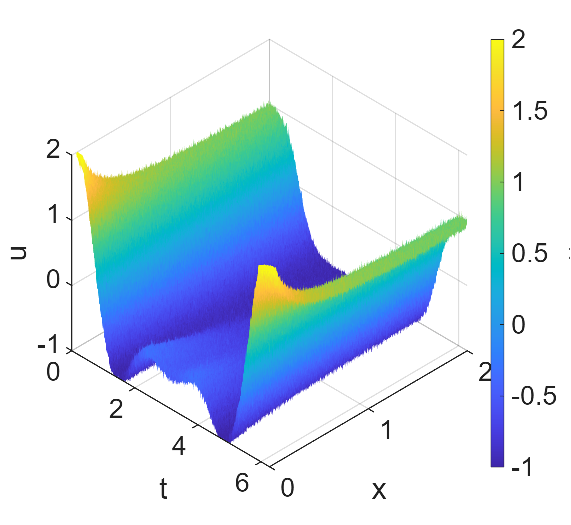}
  \label{subfig:AC_noise05_true}
\end{subfigure}\hfill
\begin{subfigure}[t]{0.32\linewidth}
  \centering

  \includegraphics[width=\linewidth]{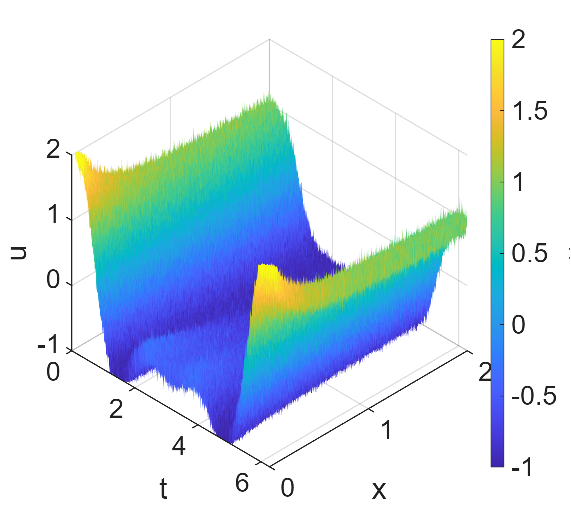}
  \label{subfig:AC_noise10_true}
\end{subfigure}

\vspace{0.5em}

\begin{subfigure}[t]{0.32\linewidth}
  \centering

  \includegraphics[width=\linewidth]{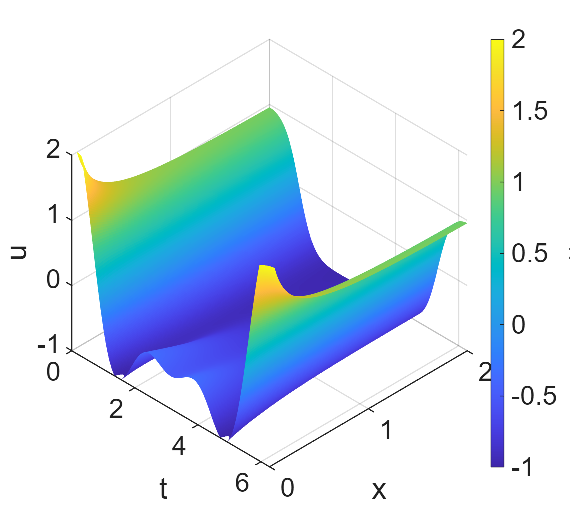}
  \label{subfig:AC_noise00_identified}
\end{subfigure}\hfill
\begin{subfigure}[t]{0.32\linewidth}
  \centering

  \includegraphics[width=\linewidth]{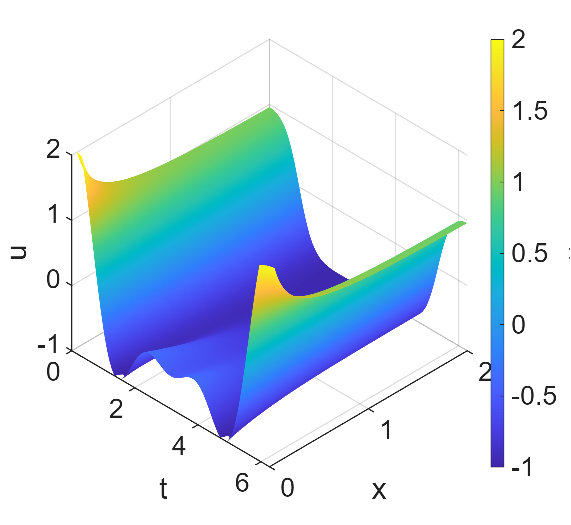}
  \label{subfig:AC_noise05_identified}
\end{subfigure}\hfill
\begin{subfigure}[t]{0.32\linewidth}
  \centering

  \includegraphics[width=\linewidth]{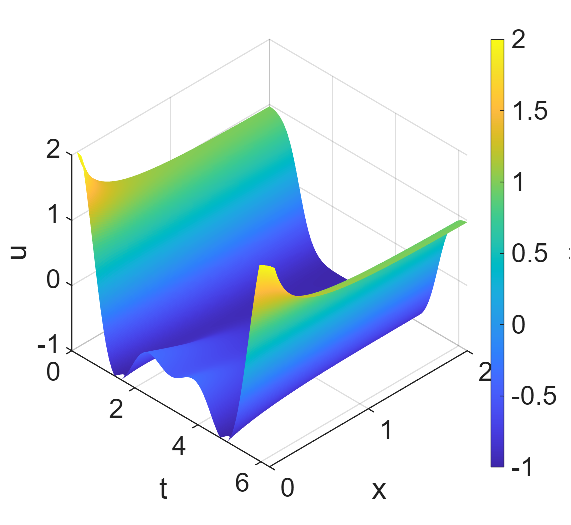}
  \label{subfig:AC_noise10_identified}
\end{subfigure}

\caption{Comparison between the true and identified trajectories of the Allen--Cahn equation under different noise levels. Top row: True trajectories. Bottom row: Identified trajectories. Columns: 0\% noise (left), 5\% noise (middle), and 10\% noise (right).
}
\label{fig:AllenCahn_2x3_grid}
\end{figure}

\noindent\textbf{Identification results and performance.}
We identify the governing PDE using both the baseline and energy-informed libraries across the four configurations. Under various noise levels, the TPR averaged over 20 trials is shown in Figure \ref{fig:AC_identification}. With the baseline library, the regression must navigate an unconstrained hypothesis space of 18 candidate terms, increasing the risk of selecting non-physical or anti-dissipative terms that violate fundamental thermodynamic principles. In contrast, the Energy--Dissipation prior restricts the search strictly to valid gradient flow dynamics via $\mathcal{D}_{\text{GF}}^{\text{valid}}$, significantly reducing the degrees of freedom. Across all noise levels, the proposed method exhibits stable performance. It correctly identifies the active support corresponding to the Allen--Cahn dynamics with an optimal sparsity of $\theta^*=3$. Notably, instead of selecting arbitrary polynomials, the method specifically recovers the gradient flow structure of the double-well potential defined in \eqref{eq:allen_cahn_energy_def}, yielding the model$$u_t = \gamma u_{xx} - \kappa (u^3 - u),$$where the identified parameters satisfy $\gamma \approx 1$ and $\kappa \approx 1$. This physically consistent recovery ensures that the data-driven model preserves the essential phase-separation mechanism, driving the system towards the stable equilibria $u = \pm 1$ while maintaining the correct energy-dissipation balance. In contrast, the unconstrained baseline approach exhibits a noticeably reduced TPR at elevated noise levels. These results demonstrate the improved performance of our proposed method: the energy prior ensures thermodynamically consistent, physically interpretable recovery, while the weak-form integration substantially enhances numerical stability and robustness to noise. Furthermore, as qualitatively illustrated in Figure \ref{fig:AllenCahn_2x3_grid} and quantitatively supported by the relative $\ell_2$ errors in Table \ref{tab:allen_cahn_E2_error}, the identified models closely reproduce the true spatiotemporal trajectories across noise levels of $\{0\%, 5\%, 10\%\}$, confirming the framework's reliability in reconstructing gradient flow structures from noisy data.

\begin{table}[ht!]
\centering
\caption{Relative $E_2$ errors of the identified Allen--Cahn solutions under different noise levels.}
\label{tab:allen_cahn_E2_error}
\renewcommand{\arraystretch}{1.2}
\begin{tabular}{c c}
\toprule
\textbf{Noise level (\%)} & \textbf{Relative $E_2$ error} \\
\midrule
0  & $6.23\times10^{-18}$ \\
5  & $3.31\times10^{-3}$  \\
10 & $1.41\times10^{-2}$  \\
\bottomrule
\end{tabular}
\end{table}

\section{Conclusions}
\label{sec:conclusions}
In this paper, we proposed a unified, prior-informed weak-form framework for data-driven PDE identification that directly embeds physical structure at the dictionary construction stage. By incorporating prior information, the candidate library is constrained to physically admissible terms, improving conditioning, interpretability, and robustness of the algorithm. Coupled with trimming and residual-reduction model selection, the pipeline yields interpretable models that are faithful to the target physics. The effectiveness of the proposed framework is demonstrated on priors including Hamiltonian (skew-gradient, energy conserving), conservation-law (flux-form with cross-directional coefficient tying), and energy–dissipation (variational, dissipative).
Extensive studies across canonical systems, such as Hamiltonian oscillators, three-body dynamics, two-dimensional shallow-water equations, and Allen–Cahn dynamics, demonstrate consistently higher true-positive rates, stable coefficient recovery, and structure-preserving behavior under substantial noise, surpassing no-prior baselines in both weak- and strong-form settings.

\section*{Acknowledgments}
We would like to thank Mr. Lingyun Deng for discussion on this project.

\bibliographystyle{abbrv}
\bibliography{references}

\end{document}